\documentclass[preprint, 3p, 11pt, number]{elsarticle}
\biboptions{sort&compress}

\usepackage{comment}
\usepackage{graphicx}
\usepackage{amssymb}

\usepackage{lineno}
\usepackage[utf8]{inputenc}
\usepackage{algorithm}
\usepackage{algorithmicx}
\usepackage{algpseudocode}
\usepackage{amsmath}
\usepackage{amsthm}

\usepackage{dsfont}
\usepackage{enumitem}
\usepackage{etaremune}
\usepackage{eurosym}
\usepackage[table]{xcolor} 
\usepackage{natbib}
\usepackage{listings}
\usepackage{bm}
\usepackage{booktabs}
\usepackage{nicefrac}
\usepackage{color, colortbl}
\definecolor{Gray}{gray}{0.9}
\usepackage{url}
\usepackage{caption}
\usepackage{subcaption}
\usepackage{microtype} 
\usepackage{multicol}
\usepackage{multirow}
\usepackage{wrapfig}
\usepackage{framed} 
\usepackage{nomencl}
\usepackage{tablefootnote}
\usepackage{arydshln}
\makenomenclature
\renewcommand*\nompreamble{\begin{multicols}{2}}
\renewcommand*\nompostamble{\end{multicols}}

\newtheoremstyle{break}
 {\topsep}{\topsep}%
 {\itshape}{}%
 {\bfseries}{}%
 {\newline}{}%
\theoremstyle{break}
\newtheorem*{remark}{Remark}

\newcommand{\beginsupplement}{%
 \setcounter{table}{0}
 \renewcommand{\thetable}{S\arabic{table}}%
 \setcounter{figure}{0}
 \renewcommand{\thefigure}{S\arabic{figure}}%
 }

\newcommand{\infolevel}[1]{\textbf{L}_{#1}}

\journal{Reliability Engineering \& System Safety}

\begin{document}

\newtheorem{example}{Example}
\newtheorem{definition}{Definition}
\newtheorem{theorem}{Theorem}
\newtheorem{proposition}{Proposition}
\begin{frontmatter}


\title{Maintenance Optimization for Asset Networks with Unknown Degradation Parameters}


\author[sor,tue]{Peter Verleijsdonk\corref{corauthor}}
\ead{p.verleijsdonk@tue.nl}
\author[ieis,tue]{Collin Drent}
\ead{c.drent@tue.nl} 
\author[sor,tue]{Stella Kapodistria}
\ead{s.kapodistria@tue.nl}
\author[ieis,tue]{Willem van Jaarsveld}
\ead{w.l.v.jaarsveld@tue.nl} 

\address[sor]{Department of Mathematics and Computer Science}
\address[ieis]{Department of Industrial Engineering and Innovation Sciences}
\address[tue]{Eindhoven University of Technology, P.O. Box 513, Eindhoven 5600 MB, the Netherlands}
\date{\today}

\cortext[corauthor]{Corresponding author}

\begin{abstract}
We consider the key challenge of maintenance optimization for asset networks where degradation parameters are \emph{heterogeneous} and \emph{unknown}, and must be inferred from real-time degradation data. Combining stochastic modeling and Bayesian statistics, we formulate a partially observable Markov decision process that incorporates estimating shock rates and sizes. This formulation allows us to analytically establish that optimal replacement policies are threshold-type in both degradation levels and parameters. Moreover, we propose an open-loop feedback approach that allows policies trained via deep reinforcement learning (DRL) with access to true parameters to remain effective when deployed using real-time Bayesian point estimates. Complementarily, we develop a Bayesian Markov decision process (BMDP) framework in which the agent maintains and updates posterior distributions during deployment, capturing parameter uncertainty evolution and enabling \emph{scalable} DRL policies that adapt as new data become available. We evaluate our approaches on synthetic data and a real-world case involving interventional X-ray filaments. The proposed DRL methods consistently outperform traditional heuristics. Policies trained for the BMDP remain robust when priors are estimated from historical data and perform well in highly heterogeneous networks. Access to true parameters provides only marginal cost improvements, highlighting the ability of the approach to make effective decisions under limited information.
\end{abstract}

\begin{keyword}
Condition-based maintenance \sep Optimization \sep Bayesian inference \sep Deep reinforcement learning 


\end{keyword}
\end{frontmatter}

\section{Introduction}\label{sec:introduction}

Optimizing the maintenance of an installed base of capital goods, such as wind turbine parks, medical imaging equipment, and wafer steppers, is critical for ensuring continuous availability and operational efficiency. The unavailability of such assets imposes significant costs, with an estimated annual economic impact of around \$50 billion~\citep{coleman, wallstreet}. Nearly half of this unavailability results from component failures, which not only cause substantial financial losses but also pose safety hazards and disrupt essential supply chains. Minimizing these failures\textemdash particularly within sectors such as energy, healthcare, and semiconductors\textemdash is therefore essential to reduce downtime costs and mitigate safety risks.

To achieve maximum availability and minimize unnecessary maintenance costs, condition-based maintenance (CBM) has become a prominent proactive maintenance strategy. In CBM, sensors continuously collect operational data and transmit them in real time. The degradation signal is constructed from these measurements to represent the health of the asset over time, enabling maintenance actions to be triggered only when failure becomes imminent. While CBM may significantly reduce unnecessary interventions and costly downtime, effectively utilizing real-time degradation data to optimally schedule maintenance across an asset network remains a non-trivial challenge.

The challenge of optimizing CBM for real-world asset networks arises from two interrelated issues: \emph{asset dependencies} and \emph{component heterogeneity}. Asset dependencies, such as economic dependencies arising from joint maintenance setup costs, shared spare parts and limited maintenance resources, significantly influence maintenance scheduling across networks of interconnected assets~\citep{OLDEKEIZER2017405}. Despite their practical significance, such dependencies are often ignored in traditional CBM optimization approaches due to the analytical complexity that results from simultaneously tracking the condition of multiple assets.

Component heterogeneity presents an additional challenge, as even nominally identical components can exhibit a significant variability in degradation behavior due to differences in manufacturing, environmental conditions, or operational usage. Traditional CBM approaches typically assume known degradation parameters and do not incorporate real-time parameter learning from observed degradation data~\citep{DEJONGE2020805, ARTS2024}. Increasingly, scholars address the challenge of real-time learning by integrating Bayesian inference into maintenance optimization to dynamically update beliefs about component degradation parameters. While promising, such Bayesian models have primarily been studied in single-asset contexts~\citep{elwany2011structured, kim2013joint, chen2015condition, van2017maintenance, msom}, where they remain analytically tractable.

Recent work in the field of maintenance optimization has made progress either in (i) integrating uncertain prognostics into maintenance planning~\citep{mitici2023dynamic, zhuang2023prognostic} and developing maintenance planning under partial observability or model uncertainty~\citep{arcieri2023bridging, tseremoglou2024condition}, or (ii) optimizing economically dependent multi-component CBM under a known degradation model~\citep{zhang2024optimal}. What remains largely missing is a unified approach in which parameter uncertainty evolves over time and is updated as new degradation observations arrive, while policies still coordinate actions across assets due to shared setup costs\textemdash such as the fixed expenses incurred when dispatching maintenance engineers to service equipment.

In this paper, we contribute: (i) a CBM model for economically dependent assets with unknown and heterogeneous degradation parameters; (ii) a dual formulation as both a partially observable Markov decision process (POMDP) and, under conjugacy assumptions, a Bayesian Markov decision process (BMDP) with a compact belief-state representation; (iii) two solution approaches based on deep reinforcement learning (DRL) for the POMDP/BMDP maintenance optimization problem: (a) open-loop deployment of policies trained under full information and (b) direct training on posterior information; (iv) analytical monotonicity results for the full-information Markov decision process (MDP) with respect to both degradation levels and parameters; and (v) numerical and case-study evidence quantifying when belief-aware training is essential and when open-loop deployment degrades under high heterogeneity.

The assets are assumed to degrade independently, with the shared setup cost constituting the only dependence between assets in the network. In addition, costs are incurred for corrective or preventive component replacements. Other factors, such as safety hazards or operational disruptions, are often incorporated indirectly through the corrective-to-preventive cost ratio~\citep{DEJONGE2020805}. Component degradation is modeled using stochastic shock models~\citep{esary1973shock, sobczyk1987stochastic}. We assign prior distributions to represent uncertainty about key degradation parameters, specifically, the shock occurrence rate and the distributional parameter for the damage incurred per shock. Conjugate priors facilitate analytically tractable Bayesian updating as degradation data accumulate, following the single-asset framework in~\citep{msom}.

To address the resulting high-dimensional control problem, we develop methodological innovations and insights that enable the training of DRL-based policies to optimize CBM for networks of interdependent assets that degrade according to a shock process with unknown parameters, which are progressively inferred from real-time degradation signals. Specifically, we develop two distinct Bayesian simulation environments for a multi-asset setting, each supporting a different DRL approach. The \emph{first modeling approach} explicitly simulates the true degradation parameters for each component; these are drawn from a population distribution at the time of component replacement. Subsequent degradation then evolves according to a stochastic shock process governed by these sampled parameters. Since the true degradation parameters are unobservable in practice, this setup constitutes a POMDP. To enable deployment of the trained policies, we propose an \emph{open-loop feedback approach}. Concretely, we collapse the posterior distributions into point estimates, which serve as input to the neural network in place of the true degradation parameters. To complement this first approach, we propose a \emph{second modeling approach} that assumes the availability of a conjugate prior and, in that case, is equivalent in \emph{objective} to the first approach\textemdash although the two differ in methodology. Instead of simulating the true degradation parameters, this approach explicitly maintains and updates uncertainty by tracking only the posterior distributions of the parameters. These posteriors are continuously updated through real-time degradation observations. This approach constitutes a BMDP, which forms the foundation for the development of DRL-based policies that operate \emph{directly on the posterior distributional information} available to the asset manager. This marks a significant advancement beyond single-asset models by jointly addressing learning and maintenance optimization in heterogeneous asset networks.

As the training algorithm for our DRL-based policies, we employ specific forms of approximate policy iteration (API)~\citep{temizoz2023deep}, in which the derived monotonicity results are incorporated both through action space constraints during training and in the design of the initialization heuristics. In particular, two heuristic policies are adopted as benchmarks: (i) a two-threshold control-limit heuristic that groups maintenance actions and is used to initialize API, and (ii) an integrated Bayes heuristic known to yield near-optimal performance in the single-asset setting.

Through extensive numerical experiments and a practical case study utilizing real-world degradation data from interventional X-ray (IXR) filaments, we demonstrate the effectiveness and practical applicability of our DRL approaches. To better understand how policy performance varies with the availability of degradation information, we extend the concept of \emph{information levels}~\citep{dtmpa} to our Bayesian setting: Information level $\infolevel{2}$ corresponds to full knowledge of the degradation parameters for each component. In information level $\infolevel{1}$, this parameter information is unknown, but the distribution of degradation parameters is available. Finally, information level $\infolevel{0}$\textemdash used only in the case study\textemdash represents a setting where even the distributional information must be estimated from historical data. Although the proposed framework is in this paper directly applied to medical equipment\textemdash specifically the filament component of the IXR system\textemdash it has broad applicability across other domains, including wind farms, wafer steppers, manufacturing processes, and beyond. The integration of Bayesian parameter learning, economic dependence, and DRL-based policy optimization applies to asset networks where degradation parameters are unknown and need to be learned over time.

Our numerical results show that DRL-based policies trained under both the POMDP and BMDP frameworks significantly outperform heuristic benchmarks when degradation parameters must be estimated (information level $\infolevel{1}$). While the POMDP-based DRL approach performs well in controlled settings with limited component heterogeneity, the BMDP-based DRL approach generally yields superior performance\textemdash particularly in scenarios with substantial heterogeneity in degradation parameters, as observed in the IXR case study. This advantage arises from the BMDP’s ability to maintain and update the belief state in real time. By leveraging the full distributional information, policies trained under the BMDP framework enable more informed and effective maintenance decisions. Furthermore, our experiments reveal that having access to the true degradation parameters ($\infolevel{2}$) yields only modest cost reductions compared to when these parameters must be inferred ($\infolevel{1}$). Finally, our case study establishes the effectiveness of our methods when distributional information on degradation parameters must be estimated from historical data. In summary, our experiments demonstrate that integrating Bayesian inference and DRL yields effective policies for CBM of asset networks, even when degradation parameters are unknown and must be inferred from real-time degradation data.

\begin{figure}[ht!]
 \centering
\includegraphics[width=1\linewidth]{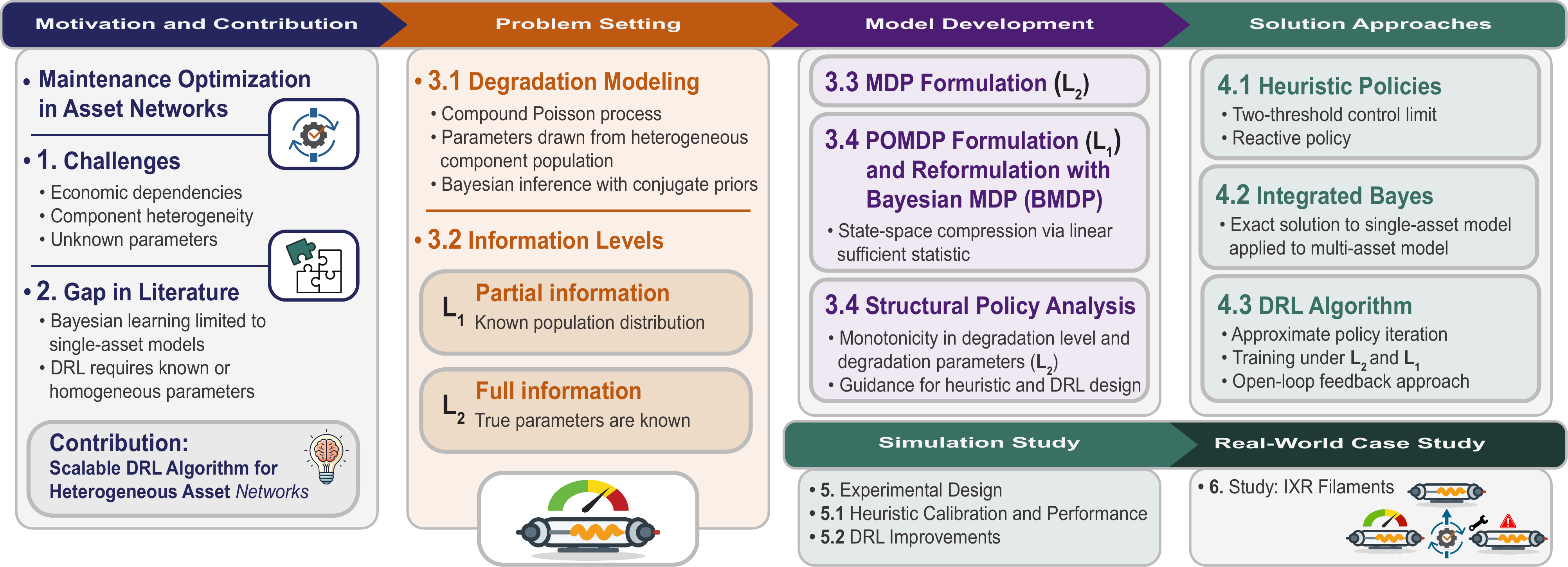}
 \caption{Flowchart illustrating the conceptual structure of the paper and the relationships between the modeling framework, theoretical results, algorithmic design, and numerical results.} \label{fig:flowchart}
\end{figure}

The remainder of the paper is organized as follows. \mbox{Section \ref{sec:literature}} presents an overview of the relevant literature. We provide a detailed model formulation in \mbox{Section \ref{sec:model}}. In \mbox{Section \ref{sec:heuristics}}, we detail the heuristic solutions and provide a brief overview of the DRL algorithm. We demonstrate the effectiveness of the proposed solutions through a concise simulation study in \mbox{Section \ref{sec:sim_study}}. In \mbox{Section \ref{sec:case_study}}, we present the results of the case study on degrading IXR filaments. We offer concluding remarks in \mbox{Section \ref{sec:conclusion}}. A visual roadmap of the paper’s structure and the relationships between its main components is provided in \mbox{Figure \ref{fig:flowchart}}.

\section{Literature review}\label{sec:literature}

Our work focuses on the optimization of CBM for asset networks with economic dependencies, the learning of maintenance strategies under component heterogeneity, and the application of DRL to CBM. We next review these three research streams.

\subsection*{Maintenance of asset networks amid economic dependencies} 

Maintenance models may be classified as single-asset or multi-asset~\citep{OLDEKEIZER2017405,DEJONGE2020805}. Multi-asset models extend single-asset models by incorporating joint maintenance policies, accounting for various dependencies such as economic, structural, stochastic, or resource dependencies. Economic dependence refers to the cost relationship between joint maintenance of multiple assets and the maintenance of individual assets. If joint maintenance is more expensive, it indicates negative economic dependence; if it is less expensive, it indicates positive economic dependence. 

Various models for CBM in systems with economic dependence have been developed in recent literature due to their relevance in industrial settings. \Citet{WIJNMALEN199752} propose a heuristic algorithm that minimizes long-run average costs in a multi-component system by determining simple, component-wise control limits and leveraging coordinated repairs to reduce shared setup costs. \Citet{BOUVARD2011601} propose a maintenance optimization method that dynamically groups tasks for multi-component systems using updated degradation data, demonstrated on commercial heavy vehicles. \citet{TIAN2011renewenerg} propose a CBM strategy for wind farms that leverages neural networks to predict the degradation level and corresponding remaining useful lifetime (RUL) from real-time data. \citet{TIAN2011ress} develop a CBM policy for multi-component systems using a proportional hazards model and a numerical algorithm to evaluate the cost of the proposed policy.

\citet{zhu2015condition} introduce a novel CBM policy for multi-component systems experiencing continuous stochastic degradation, utilizing a joint maintenance interval to reduce shared setup costs. \citet{OLDEKEIZER2016531} develop a dynamic programming model to determine the optimal CBM strategy for systems with both economic dependencies and redundancy, demonstrating that it significantly outperforms heuristic policies and offers insights into the optimal policy structure. \citet{OLDEKEIZER2018319} investigate the trade-off between prompt replacement of failed components and maintenance clustering opportunities in redundant systems with economic dependencies, revealing that heuristic threshold policies may be suboptimal and that misinterpreting load-sharing effects can lead to higher maintenance expenses. \citet{DO201986} present a model for a CBM policy in a two-component system with stochastic dependencies, where the degradation of each component is influenced by the other's state, and economic dependencies, where combined maintenance activities are shown to be more cost-effective. Similarly,~\citet{OAKLEY2022108321} propose a CBM policy for multi-component systems that balances the trade-offs between the urgency of replacing failed components due to increased load and the economic benefits of clustering replacements to minimize downtime and maintenance costs. \Citet{zhang2024optimal} incorporates imperfect maintenance into CBM optimization for multi-component systems using genetic algorithms. We refer the reader to~\citet[\mbox{Table 1}]{zhu2015condition} for a more extensive summary of CBM models proposed for multi-component systems.

The geographical layout of asset networks often introduces a complex dependency. \Citet{abdul2018optimally} investigate the challenge of optimally replacing multiple stochastically degrading systems within a shared environment using CBM. \Citet{soltani2023structured} study the problem of optimally maintaining an offshore wind turbine farm under both economic and stochastic dependence due to shared maintenance setup costs and their common environment. They establish monotonicity of the cost function jointly in the degradation level and environmental state, characterize the structure of the optimal replacement policy, and show that sharing maintenance resources is cost-effective. \Citet{LEPPINEN2025103162} introduce directed graphs to represent the economic and structural dependencies of a multi-component system, including scenarios where maintenance on one component may require the disassembly or maintenance of others, and solve the resulting MDP using a modified policy iteration algorithm to determine the most cost-efficient policy.

\subsection*{Learning maintenance strategies under component heterogeneity} 

Component heterogeneity, characterized by varying degradation rates and failure patterns, necessitates the integration of learning the component-specific characteristics of the degradation process and optimizing maintenance strategies. Degradation processes can be represented using either a discrete or continuous set of degradation states, with Markov chains commonly used to model the degradation process in the discrete case. Information gathered from sensors can enhance the effectiveness of scheduled inspections. 

Most existing studies on uncertainty in degradation processes typically assume a specific parametric form, with uncertainty modeled through its parameters and captured via a prior distribution. Bayesian inference is then used to update this uncertainty over time\textemdash but so far, these approaches have been limited to single-asset systems only. \citet{van2017maintenance} develop a POMDP model to incorporate population heterogeneity for maintenance scheduling of a single-asset system that stochastically degrades; however, the population of spare parts consists of multiple indistinguishable types that degrade at varying rates. \citet{elwany2011structured} and~\citet{SiWiener} both study a Wiener degradation process with an unknown drift parameter. In both studies, periodic inspections are used to estimate the parameter, with Bayesian inference applied in the former and maximum likelihood estimation in the latter. These estimation methods, combined with observed degradation levels, guide the decision-making process for maintenance interventions. \citet{FLAGE201216} consider a single asset subject to a degradation process with unknown parameters. Sequential inspections assess the degradation level, with the timing of the next inspection based on the current degradation level. Uncertainty in the parameters is modeled using a prior distribution and updated in a Bayesian manner. \Citet{msom} investigate a single-asset model with one critical component that degrades according to a compound Poisson process with unknown parameters, leveraging real-time degradation data to infer the component’s degradation behavior and adjust decision-making accordingly. They demonstrate that, accounting for component heterogeneity, the optimal policy depends on both the asset’s age and the observed degradation signal, and that integrating learning with decision-making yields significant cost reductions.

In contrast to \Citet{msom} and other single-asset models described above, our work considers a network of economically dependent assets, where managers must account for shared setup costs and coordinate maintenance decisions across multiple assets. Extending such single-asset models to this multi-asset setting substantially increases computational complexity, rendering exact solution approaches for the resulting BMDP intractable. We therefore develop scalable DRL-based policies that enable integrated learning and maintenance optimization in heterogeneous asset networks under parameter uncertainty.

\subsection*{Deep reinforcement learning for condition-based maintenance}

DRL-based approaches to maintenance optimization have demonstrated promising results across various settings. \Citet{kuhnle2019reinforcement} train opportunistic DRL-based maintenance strategies using proximal policy optimization for parallel assets, achieving reductions in downtime and costs compared to traditional strategies such as reactive and time-based maintenance. \Citet{ZHANG2020107094} propose a DRL algorithm based on deep Q-network (DQN) to train policies for CBM in multi-component systems with stochastic and economic dependencies, and show in a numerical study that the trained policies outperform benchmarks set by heuristic policies. \Citet{MOHAMMADI2022108615} develop a DRL-based method using double DQN (DDQN) for the joint optimization of maintenance and renewal planning under practical constraints. Using historical inspection and maintenance data to simulate a rail infrastructure environment, they show that the approach generates cost-effective policies that enhance network reliability and safety. \citet{Hung2024} employ DDQN to train DRL-based maintenance policies in a stochastic factory setting characterized by various degradation levels, uncertain repair times and fluctuating machine workloads. 

\Citet{LEE2023108908} develop an integrated predictive maintenance framework that combines RUL estimation via convolutional neural networks with a DRL-based maintenance policy trained using a soft actor–critic algorithm. The policy triggers maintenance actions based on the estimated RUL distribution and yields significant cost savings and downtime reduction for aircraft turbofan engines. \Citet{tseremoglou2024condition} propose a two-stage CBM framework for aircraft fleets, where RUL predictions are first used to construct a maintenance policy via a POMDP and then integrated into a rolling-horizon DQN to schedule preventive and corrective tasks subject to resource constraints. \Citet{zhuang2023prognostic} proposes a prognostics-driven framework that uses Bayesian deep learning to generate uncertainty-aware RUL distributions and dynamically updates maintenance and spare-part ordering decisions for turbofan engines.

\citet{dtmpa} propose a DRL approach based on distributional DQN to minimize maintenance and downtime costs in asset networks serviced by a single maintenance engineer, where degradation levels are only partially observable. Building upon this work,~\citet{kdtmpa} develop DRL-based solutions trained via a form of API\textemdash shown to scale better than the aforementioned training algorithms\textemdash for industrial asset networks maintained by multiple engineers; however, their approach assumes full observability of all degradation levels.

Recent studies have made important progress (see, e.g.,~\citet{andriotis2020deep}) by integrating Bayesian inference with DRL in the context of maintenance optimization. These efforts focus on inferring the degradation state of components rather than the parameters of the underlying degradation process. Thus, while such studies represent key steps toward maintenance optimization under uncertainty, they typically assume that the parameters of the degradation processes are known and can be used to define the simulation environment in which the DRL-based policy is trained.

Our work represents a next step in the integration of learning and maintenance optimization by addressing degradation parameter uncertainty in multi-asset systems. We develop a modeling framework in which the degradation process itself is partially unknown and must be learned in real time from degradation signals. This enables DRL-based policies to adapt not just to the observed condition of assets, but also to the underlying heterogeneity in their degradation behavior. While such models have been considered in single-asset systems, we appear to be the first to consider this learning problem in the realistic context of (economically) interdependent asset networks, where coordination is essential due to, e.g., shared maintenance setup costs. We formalize the resulting optimization problem as both a POMDP and a BMDP, and propose scalable solutions that operate on evolving posterior distributions of degradation parameters. This allows us to make a significant next step in data-adaptive maintenance policies for realistic, industrial-scale asset networks.
\section{Model formulation}\label{sec:model}

In this section, we formulate the maintenance problem for asset networks, accounting for economic dependencies and component heterogeneity.

\subsection{Compound Poisson degradation}

We consider a set of assets (machines) $\mathcal{M} = \{1, \ldots, M\}$, each with a single critical component that degrades independently as random shocks arrive. From this point forward, we will refer to the degradation of components as asset degradation. Shocks occur according to a Poisson process, and the damage that accumulates during each shock is a non-negative random variable, meaning that degradation follows a compound Poisson process. This type of shock model is appropriate for, e.g., certain metal and ceramic components in trains, aircraft, and medical equipment (including the IXR systems of our case study) that primarily deteriorate when subjected to discrete stress events~\citep{esary1973shock}. In such settings, degradation accumulates through a sequence of shocks associated with random events, such as thermal cycles, mechanical loads, or electrical stress. Modeling these shock arrivals as a Poisson process is common in the reliability and maintenance literature, as it provides a tractable representation of random and approximately independent stress events occurring over time. Indeed, shock-based degradation models with Poisson arrivals have been widely used in maintenance optimization; see, e.g.,~\citep{tian2021optimal, msom, drent2024optimal, drent2024condition, wang2024degradation}. Consistent with this literature, we adopt the same modeling approach. The Poisson intensity of shock arrivals at machine $m \in \mathcal{M}$ is denoted by $\lambda_m \in \mathbb{R}_+$. The random damage amount at machine $m$ resulting from a shock adheres to the distribution of a member of the one-parameter (denoted as $\phi_m$) exponential family. The probability density or mass function of such a random variable can be written as
\begin{equation*}
 f_m(x \mid \phi_m) = h_m(x)e^{\phi_m T_m(x)-A_m(\phi_m)}, 
\end{equation*}
where $T_m(x)$ represents the sufficient statistic, and $h_m(x)$ and $A_m(\phi_m)$ are known functions. We assume that shock sizes are non-negative and that the sufficient statistic is linear, i.e., $T_m(x) \equiv x$, which enables a state space reduction in our optimization problem. The function $T_m(\cdot)$ contains all the information necessary to compute any estimate of the parameter $\phi_m$. In the literature, this group of distributions is commonly known as the linear exponential family or natural exponential family, named for its linear sufficient statistic, and was first introduced by~\citet{morris}. To illustrate our approach, we consider the geometric distribution (supported on $\mathbb{Z}_+$) as a representative example. 

\begin{example}[Geometric distribution]\label{example:geo}
 The probability mass function of a geometrically distributed shock size (with support $\mathbb{Z}_+$) with parameter $p_m \in (0,1)$ is given by
\begin{equation*}
 f(x \mid p_m ) = (1-p_m)^x p_m = e^{ \ln(1-p_m)x - \ln(1/p_m)}.
\end{equation*}
Note that $h_m(x) = 1$, $T_m(x) = x$, $\phi_m = \ln(1-p_m)$ and $A_m(\phi_m) = \ln(1/p_m)$.
\end{example}

\begin{remark}
The geometric distribution serves as an illustrative example, but the framework extends to any one-parameter exponential family with linear sufficient statistic $T_m(x)=x$. This class includes several commonly used discrete and continuous distributions for modeling degradation increments. 
In the discrete case, examples of the linear exponential family include the Poisson distribution, the binomial distribution (with fixed number of trials), and the negative binomial distribution (with fixed number of successes). 
In the continuous case, the exponential distribution, normal distribution (with fixed variance), inverse Gaussian (with fixed shape parameter), and the gamma distribution (with fixed shape parameter) are all members of the linear exponential family. In \ref{app:examples}, we demonstrate for two standard discrete distributions\textemdash the binomial and the Poisson\textemdash how they can be written in canonical exponential-family form with linear sufficient statistic, thereby fitting within our modeling framework.
\end{remark}

We assume continuous, perfect, remote access to the degradation level for each machine through sensory equipment, while interaction with the system is limited to evenly spaced decision epochs (without loss of generality, we rescale time such that the time between two consecutive decision epochs equals one) corresponding to scheduled maintenance opportunities. Two remarks are in order to motivate this assumption. First, regular maintenance intervals are common in practice when maintenance planning and logistics must be arranged well in advance. For example, in geographically dispersed asset networks, maintenance visits are often scheduled at fixed intervals (e.g., semi-annually), as mobilizing crews and spare parts requires substantial preparation. Second, with the increasing deployment of internet-connected assets, embedded sensors, and digital monitoring infrastructures, near-continuous and accurate condition monitoring is becoming standard in many industrial applications. We therefore deliberately study this data-rich setting as a benchmark for integrated learning and decision-making. In the case study presented in \mbox{Section \ref{sec:case_study}}, a time unit is defined to roughly correspond to the minimum time required to dispatch a maintenance team or deliver a spare part.

Let $N_m(t)$ denote the total number of shocks incurred by the $m$-th asset since its installation up to its \emph{operational age} (that is, the time since the last replacement) $t\in\mathbb{R}_+$, i.e., $N_m(t)$ is a Poisson process with rate $\lambda_m$. Here, it is important to note that real-time continuous monitoring of the machine, rather than relying solely on periodic inspections, enables an accurate count of the total number of shocks. The $m$-th asset's degradation level is denoted by $X_m(t) \in [0, \xi_m]$ for some failure threshold $\xi_m > 0$. When the asset's degradation level reaches or exceeds $\xi_m$, it breaks down and requires corrective maintenance at the next decision epoch. Maintenance is instantaneous and restores the asset to an as-good-as-new condition.

We denote the number of shocks that arrive in the time interval $(t-1,t]$ by $K_m\big((t-1,t]\big) = N_m(t) - N_m(t-1)$. Furthermore, let $Y^{(i)}_m$ denote the size of the $i$-th shock at machine $m$ since the last replacement of the asset. The total incurred damage $X_m(t)$ is a compound Poisson process and satisfies 
\begin{equation*}
 X_m(t) = \sum_{i=1}^{N_m(t)} Y^{(i)}_m,
\end{equation*}
where $X_m(0) = 0$ and $N_m(0) = 0$ by definition. Furthermore, let $\bm{Y}_m\big( (t-1,t] \big) = (Y^{(N_m(t-1)+1)}_m, \ldots, Y^{(N_m(t))}_m)$ be the sizes of the shocks that arrive between age $t-1$ and $t$, and let $Z_m\big( (t-1, t] \big) = \sum_{i=N_m(t-1)+1}^{N_m(t)} Y^{(i)}_m $ be the corresponding total incurred damage. Let $k_m(t)$ denote the observed number of shocks of the $m$-th asset, that is, $k_m(t)$ is the realization of $K_m\big( (t-1,t] \big)$. Denote with $\bm{y}_m(t) = (y_m^{(1)}, \ldots, y_m^{(k_m(t))})$ the corresponding observed shock sizes. That is, $\bm{y}_m(t)$ is the realization of $\bm{Y}_m\big( (t-1,t] \big)$. The tuple $\bm{\theta}_m(t) = (k_m(t), \bm{y}_m(t))$ is the observed degradation signal of machine $m$ between ages $t-1$ and $t$. Given the assumption of a linear sufficient statistic, we do not need to keep track of individual shock sizes. Instead, we can collapse the state space by summarizing the signal as $(k_m(t), z_m(t))$, where $z_m(t) = \sum_{i=1}^{k_m(t)} y_m^{(i)}(t)$ represents the accumulated damage between ages $t-1$ and $t$. That is, $z_m(t)$ is the realization of $Z_m\big( (t-1, t] \big)$. This state space collapse is a key simplification that enables more tractable inference.

Recall that each component stems from a distinct heterogeneous population that consists of components with different degradation parameters $\lambda_m$ and $\phi_m$. In practice, these parameters are hidden; only the observed degradation signal $\bm{\theta}_m(t)$ is available at operational age $t$. In this case, we will employ a POMDP to model the integrated challenge of learning the degradation parameters while determining the optimal timing for replacement. In the \emph{underlying} MDP, these parameters $\lambda_m$ and $\phi_m$ are drawn from known distributions, $\Lambda_m$ and $\Phi_m$, respectively, and are observed by an oracle that is aware of the true population heterogeneity. The asset manager's knowledge about the degradation model and its parameters gives rise to various information levels, which are detailed in the next section.

\subsection{Information levels}\label{subsec:loi}

Inspired by the approach of~\citet{dtmpa}, we formalize the asset manager's knowledge of the degradation model and its underlying parameters using three distinct \emph{levels of information}:
\begin{enumerate}[label=($\textbf{L}_\arabic*$)]
\setcounter{enumi}{-1}
\item The asset manager has no information about the model's parameters $\lambda_m, \phi_m$, nor any knowledge of their distribution (i.e., the distributions $\Lambda_m$ and $\Phi_m$ are unknown) for each $m \in \mathcal{M}$. \label{L0}

\item The asset manager has no information about the model's parameters $\lambda_m,\phi_m$, but has full knowledge of their distribution (i.e., the distributions $\Lambda_m$ and $\Phi_m$ are known) for each $m \in \mathcal{M}$. \label{L1}

\item The asset manager has full information about the model's parameters $\lambda_m,\phi_m$ for each $m \in \mathcal{M}$. \label{L2}
 \end{enumerate}

Given an information level $\infolevel{} \in \{ \infolevel{0}, \infolevel{1}, \infolevel{2} \}$, the objective is to devise a policy $\pi^\infolevel{}$ that minimizes the total expected discounted cost of managing a specific network of assets. The baseline information level, $\infolevel{0}$, reflects the typical conditions encountered by the asset manager in practice. To address the heterogeneity of the component population, hyperparameters of the distributions $\Lambda_m$ and $\Phi_m$ for each $m \in \mathcal{M}$ should be estimated from available historical degradation data. Subsequently, the asset manager formulates and solves an approximate problem under the enhanced information level $\infolevel{1}$, and implements the resulting solution in practice. In contrast, information level \(\infolevel{2}\) provides the most accurate and informed basis for decision-making. Consequently, an optimal policy derived under this level of information achieves the lowest total expected discounted cost. This progression from $\infolevel{0}$ to $\infolevel{2}$ highlights the trade-off between the costs of acquiring additional information and the benefits of improved CBM planning. However, although relevant, these acquisition costs are not explicitly incorporated in our model and are left for future work.

The forthcoming sections formalize the model under information level $\infolevel{2}$ as an MDP and under information level $\infolevel{1}$ as a POMDP or a BMDP pending additional assumptions.

\subsection{Partially observable Markov decision process formulation}\label{subsec:pomdpmodel}

Under information level $\infolevel{2}$, a state $h$ of the asset network can be represented by a vector $h=(x_1, \lambda_1, \phi_1, \ldots, x_M, \lambda_M, \phi_M)$, with a minor abuse of notation. Here, $x_m$ denotes the degradation level of the $m$-th asset, and $\lambda_m$ and $\phi_m$ are its current degradation parameters. At every decision epoch, the asset manager must choose for each asset whether to maintain the asset or to postpone maintenance activities. Maintenance on failed assets is mandatory. Hence, the \textit{state-dependent action set for the $m$-th asset} is
$$\mathcal{U}_m(h) = \begin{cases}
\{0,1\} & \textrm{if } x_m < \xi_m,\\
\{1\} & \textrm{if } x_m \geq \xi_m.
\end{cases}$$
The \emph{state-dependent action set} $\mathcal{U}(h) = \mathcal{U}_1(h) \times \ldots \times \mathcal{U}_M(h)$ is formed by taking the Cartesian product of the $M$ individual state-dependent action sets. For tractability, we assume that maintenance capacity is sufficient to service all assets simultaneously.

The asset manager incurs costs related to either the corrective or preventive replacement of an asset. If the degradation level at a decision epoch is less than $\xi_m$, then we can either maintain the asset preventively at cost $c_m^\textrm{PM}$ or proceed to the next decision epoch without incurring any cost. If the degradation level of machine $m$ at a decision epoch is greater than or equal to the failure threshold $\xi_m$, then the failed asset is replaced correctively at cost $c_m^\textrm{CM}$. We assume that for all $m\in\mathcal{M}$ it holds that $0 < c_m^\textrm{PM} < c_m^\textrm{CM} < \infty$ to avoid unrealistic cases. This cost structure is commonly adopted in the maintenance literature~\citep{DEJONGE2020805}. Moreover, both types of replacements take negligible time, which is a reasonable assumption because in practice, replacement times are relatively small compared to the time between decision epochs. The assets are economically coupled through a shared setup cost, a widely adopted assumption for modeling economic dependence in multi-component maintenance systems~\citep{OLDEKEIZER2017405}. In practice, such setup costs represent fixed expenses associated with initiating maintenance activities, such as dispatching a maintenance crew, preparing equipment, or temporarily shutting down part of the system. Specifically, if \emph{at least} one asset is replaced during a decision epoch, a one-time setup cost $c^\textrm{ST} \geq 0$ is incurred. Although we focus on this particular economic dependence, the algorithmic developments presented in this paper readily extend to other forms of economic dependencies, such as more complex setup cost functions. After an asset $m$ is replaced, new degradation parameters $\lambda_m$ and $\phi_m$ are drawn from their respective distributions, $\Lambda_m$ and $\Phi_m$. Therefore, the costs incurred when taking action $a = (a_1, \ldots, a_M) \in \mathcal{U}(h)$ in state $h$ are:
\begin{align*}
C(h,a) &= c^\textrm{ST}\mathds{1}_{\{ \sum_{m\in\mathcal{M}} a_m > 0 \}} + \sum_{m\in\mathcal{M}}\left(c_{m}^\textrm{CM}\mathds{1}_{\{a_m = 1, x_{m} \geq \xi_m\}} + c_{m}^\textrm{PM}\mathds{1}_{\{a_m = 1, x_{m} < \xi_m\}}\right).
\end{align*}
Under information level $\infolevel{2}$, the objective of the \emph{underlying} MDP is as follows: We are interested in a policy, say $\pi^{\infolevel{2}}$, which minimizes the total expected discounted cost. A policy is defined as a series of decision rules, i.e., $\pi^{\infolevel{2}} = (\pi^{\infolevel{2}}_1, \pi^{\infolevel{2}}_2, \ldots,\pi^{\infolevel{2}}_t, \ldots )$, where the decision rule $\pi^{\infolevel{2}}_t$ at time $t$ represents a probability distribution over the action set $\mathcal{U}(h(t))$ given the state $h(t)$. Let $J(\pi^{\infolevel{2}})$ denote the total expected discounted cost, given a discount factor $\gamma \in [0,1)$. The objective is to find an optimal policy $\pi^{\infolevel{2}}_*$ that satisfies
\begin{equation}\label{eq:objectiveMDP}
 \pi^{\infolevel{2}}_* = \textrm{arg }\underset{\pi^{\infolevel{2}}}{\textrm{min }} J(\pi^{\infolevel{2}}) = \textrm{arg }\underset{\pi^{\infolevel{2}}}{\textrm{min }} \lim_{T\to\infty} \mathbb{E}_{\pi^{\infolevel{2}}} \Bigg[ \sum_{t=0}^T \gamma^t C\left(h(t),a(t)\right) \Bigg| h(0)=h\Bigg],
\end{equation}
where $(h(t), a(t))$ represents the tuple of the underlying MDP state and the corresponding action selected by the policy $\pi^{\infolevel{2}}_t$ at time $t$, $t\geq0$, and $C(\cdot)$ indicates the associated costs (maintenance and setup).

Under information level $\infolevel{1}$, the objective becomes to find a policy $\pi^{\infolevel{1}}_*$ that satisfies
\begin{equation}\label{eq:objectivePOMDP}
 \pi^{\infolevel{1}}_* = \textrm{arg }\underset{\pi^{\infolevel{1}}}{\textrm{min }} J( \pi^{\infolevel{1}} ) = \textrm{arg }\underset{ \pi^{\infolevel{1}}}{\textrm{min }} \lim_{T\to\infty} \mathbb{E}_{\pi^{\infolevel{1}}} \Bigg[ \sum_{t=0}^T \gamma^t C\left(o(t),a(t)\right) \Bigg| o(0)=o\Bigg],
\end{equation}
where $(o(t),a(t))$ denotes the tuple of the POMDP state and the corresponding action selected by the policy $\pi^{\infolevel{1}}_t$ at time $t$. The state $o(t) = (t_1, \Theta_1((t-t_1, t]), \ldots, t_M, \Theta_M((t-t_M, t]))$ contains the machine ages and observed degradation signal history for all assets, i.e., $\Theta_m((t-t_m, t]) = \{ \bm{\theta}_m(\tau) \mid \tau \in (t-t_m, t]\}$. A belief distribution over the unobserved parameters $(\lambda_1, \phi_1, \ldots, \lambda_M, \phi_M)$ must be computed from the available history.

A Markovian belief state enables a POMDP to be formulated as an MDP, where each belief represents a state. Although computing belief updates is generally computationally intractable, conjugate pairs ensure that posterior belief states remain within the same distributional family as the prior. This consistency streamlines belief updates and enables us to reformulate the POMDP as a BMDP, improving the tractability of solving the objective in \mbox{Eq. (\ref{eq:objectivePOMDP})}.

\subsection{Bayesian Markov decision process formulation}\label{subsec:bmdp}

BMDPs are used in decision-making when the true state of the environment is hidden. Instead of making decisions based on the actual state, which is hidden, managers base decisions on a belief state. A belief state is a probability distribution over all possible states that represents the knowledge about the actual state. This framework combines MDPs and Bayesian inference by updating belief states with new information retrieved from observations and actions.

The use of conjugate pairs refines this integration by ensuring that the updated belief distributions remain within the same family as the prior belief distributions. For example, the gamma distribution acts as a conjugate prior for the (unknown) rate parameter of a Poisson distribution. Thus, for each asset $m \in \mathcal{M}$, we assume that $\lambda_m$ is drawn from a gamma distribution, $\Lambda_m \sim \textrm{Gamma}(\alpha_m, \beta_m)$, where $\alpha_m > 0$ is the shape parameter and $\beta_m > 0$ is the scale parameter. Additionally, for each $m \in \mathcal{M}$, we assume that $\Phi_m$ is distributed according to the general prior for a member of the exponential family characterized by hyperparameters $r_m > 0$ and $s_m > 0$. A member of the exponential family has a conjugate prior with a density that can be expressed as
\begin{equation*}
 f_{\Phi_m}\left(\phi_m \mid r_m(t), s_m(t)\right) = H_m\left(r_m(t), s_m(t)\right)e^{r_m(t) \phi_m - s_m(t) A_m(\phi_m)},
\end{equation*}
where $H_m\left(r_m(t), s_m(t)\right)$ is a normalizing constant~\citep{ghosh2007troduction}. The beta distribution serves as a conjugate prior for the parameter of several distributions in the exponential family, for instance the geometric distribution.

When asset $m$ is installed, the parameters $\lambda_m$ and $\phi_m$ of the compound Poisson degradation process are drawn from the \emph{known} distributions $\Lambda_m$ and $\Phi_m$, and there is no available history. Thus, we interpret the parameters $\lambda_m$ and $\phi_m$ as random variables, denoted by $\tilde{\Lambda}_m$ and $\tilde{\Phi}_m$.

\paragraph{Initial belief:} At time $t=0$, the manager has a prior belief about the degradation process parameters of each asset, which is expressed as a probability distribution over all possible combinations. The joint prior belief distribution of asset $m$ satisfies 
\begin{equation*}
f^{(0)}_{\tilde{\Lambda}_m, \tilde{\Phi}_m}(\lambda_m, \phi_m) := f_{\Lambda_m}(\lambda_m \mid \alpha_m, \beta_m )\cdot f_{\Phi_m}(\phi_m \mid r_m, s_m).
\end{equation*}

\paragraph{Update with observations and actions:} At time $t$, we use the observed degradation signal $\bm{\theta}_m(t)$ to infer the joint distribution of $\tilde{\Lambda}_m$ and $\tilde{\Phi}_m$. As shown by~\citet[Proposition 1]{msom}, this joint distribution can be factorized into two independent distributions of the same form, with parameters updated solely based on the information contained in $\bm{\theta}_m(t)$ observed in the previous time period.

\begin{proposition}\label{prop:bayes_update}
The joint posterior distribution at time $t$ of $\tilde{\Lambda}_m$ and $\tilde{\Phi}_m$ is given by
\begin{equation}\label{eq:prior_distr}
 f^{(t)}_{\tilde{\Lambda}_m, \tilde{\Phi}_m}(\lambda_m, \phi_m) = f^{(t-1)}_{\tilde{\Lambda}_m}(\lambda_m \mid \alpha_m + k_m(t), \beta_m + t_m(t) )\cdot f^{(t-1)}_{\tilde{\Phi}_m}(\phi_m \mid r_m + x_m(t), s_m + k_m(t) ).
\end{equation}
\end{proposition}

This is an iterative updating process. After each action and observation, the belief distribution for each asset is updated, and these updated beliefs define the prior for the next step. Thus, at time $t$, the prior belief distribution of $\lambda_m$ is a gamma distribution with parameters 
\begin{equation}\label{eq:alpha_update}
\alpha_m(t) = \alpha_m + k_m(t)
\end{equation}
and
\begin{equation}\label{eq:beta_update}
\beta_m(t) = \beta_m + t_m(t).
\end{equation}
Similarly, the prior belief distribution of $\phi_m$ at time $t$ is a beta distribution with parameters 
\begin{equation}\label{eq:a_update}
r_m(t) = r_m + z_m(t)
\end{equation}
and
\begin{equation}\label{eq:b_update}
s_m(t) = s_m + k_m(t).
\end{equation}
Thus, the state at time $t$ of the associated BMDP can be represented by a vector $\tilde{h}(t)=(x_1(t), k_1(t), t_1(t), \ldots, x_M(t), k_M(t), t_M(t))$. The objective of \mbox{Eq. (\ref{eq:objectivePOMDP})} is equivalent to 
\begin{equation}\label{eq:objectiveBMDP}
 \pi^{\infolevel{1}}_* = \textrm{arg }\underset{\pi^{\infolevel{1}}}{\textrm{min }} J( \pi^{\infolevel{1}} ) = \textrm{arg }\underset{ \pi^{\infolevel{1}}}{\textrm{min }} \lim_{T\to\infty} \mathbb{E}_{\pi^{\infolevel{1}}} \Bigg[ \sum_{t=0}^T \gamma^t C\left(\tilde{h}(t),a(t)\right) \Bigg| \tilde{h}(0)=\tilde{h}\Bigg],
\end{equation}
where $(\tilde{h}(t),a(t))$ denotes the tuple of the BMDP state and the corresponding action selected by the policy $\pi^{\infolevel{1}}_t$ at time $t$.

\subsection{Structural properties of the underlying Markov decision process}\label{subsec:structural_properties}

In this section, we establish several structural properties of optimal replacement policies $\pi^{\infolevel{2}}_*$ of the \emph{underlying} MDP of the POMDP model presented in \mbox{Section \ref{subsec:pomdpmodel}}. We provide the proof of \mbox{Theorem \ref{thm:deglvls}} and \mbox{Theorem \ref{thm:degparams}} explicitly for the case $M=2$ in \ref{app:proofs}. The argument extends analogously to the case $M>2$.

\begin{theorem}[Monotonicity in degradation levels]\label{thm:deglvls}
Let $x = (x_1, \lambda_1, \phi_1, \ldots, x_M, \lambda_M, \phi_M)$ denote a state. Define $\mathcal{M}(x) \subseteq \mathcal{M}$ to be the set of machines for which it is optimal to do maintenance in state $x$. For any $y = (y_1, \lambda_1, \phi_1, \ldots, y_M, \lambda_M, \phi_M)$ that represents a state with more severe degradation than $x$, i.e.,
\begin{enumerate}
 \item $y_m \geq x_m$ for all $m \in \mathcal{M}(x)$, and
 \item $y_m = x_m$ otherwise.
\end{enumerate}
Then, it holds that
\begin{equation*}
 a^*(x) \in \textrm{arg} \underset{a\in\mathcal{U}(y)}{\textrm{min}} Q(y,a).
\end{equation*}

Here, $a^*(x) \in \mathcal{U}(x)$ and $a^*(y) \in \mathcal{U}(y)$ denote the optimal actions in states $x$ and $y$, respectively, and $Q(\cdot,\cdot)$ denotes the state-action value function of an optimal policy. In other words, the optimal action $a^*(x)$ for state $x$ is also an optimal action for a state $y$ with more severe degradation.
\end{theorem}

\begin{figure}[!ht]
 \centering
 \includegraphics[width=\linewidth]{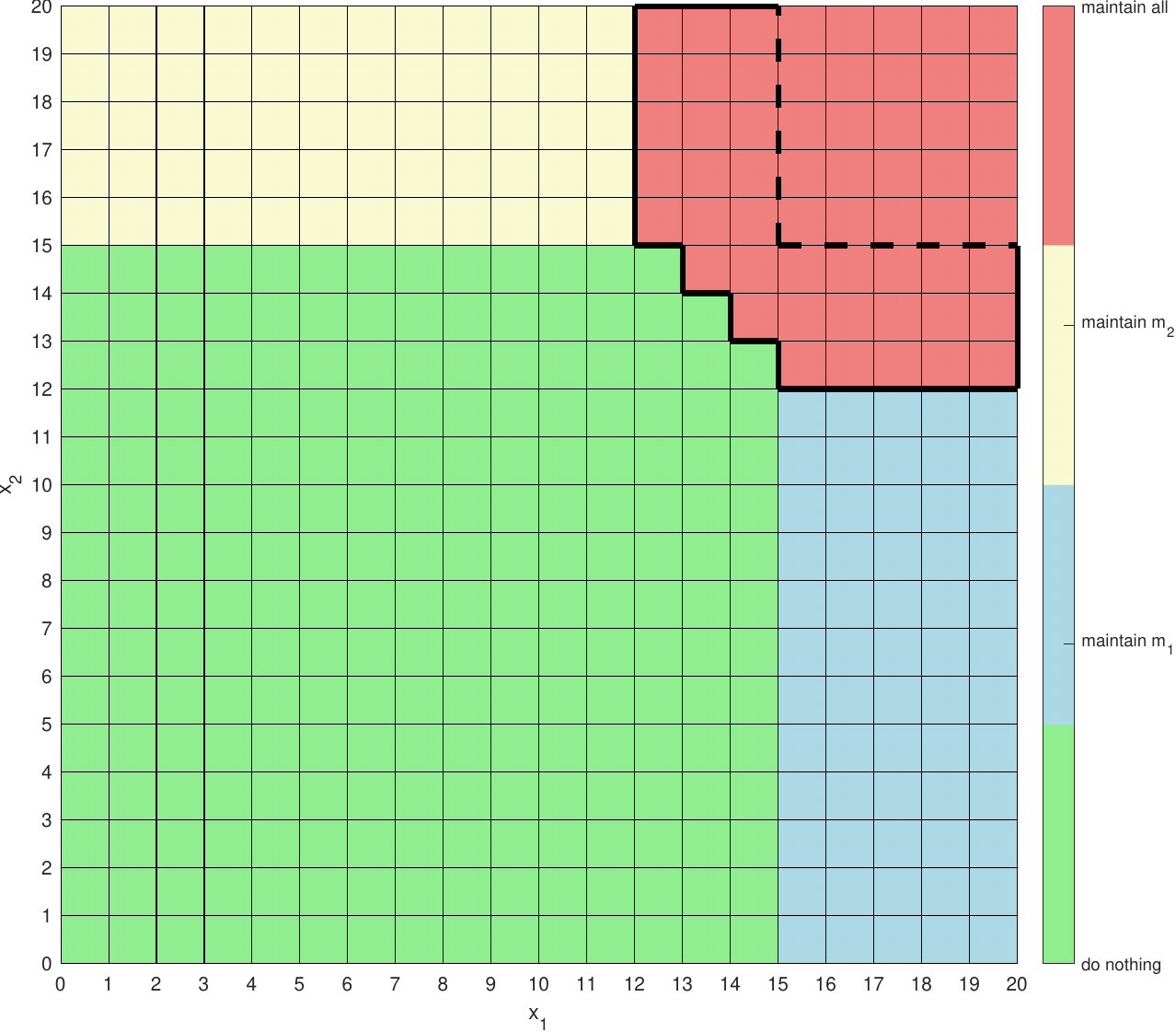}
 \caption{An optimal policy $\pi^{\infolevel{2}}_*$ of the underlying MDP for the 2-asset instance with parameters $(c_m^\textrm{PM}, c_m^\textrm{CM}, c_m^\textrm{ST}, \gamma) = (1, 10, 3, 0.99)$, $\Lambda_m \overset{\textrm{d}}{=} 5$ and $Y^{(i)}_m \overset{\textrm{d}}{=} 1$.}
 \label{fig:opt_policy}
\end{figure}

From \mbox{Theorem \ref{thm:deglvls}}, it follows that optimal policies are \emph{state-dependent} threshold policies. \mbox{Figure \ref{fig:opt_policy}} shows an example of an optimal policy for a 2-asset instance retrieved via policy iteration. Note that, due to the presence of shared setup costs, the optimal action for all states within the black boundary is to ``maintain all'', despite the fact that the optimal PM threshold for an asset considered in isolation is $15$. The solid black boundary is inclusive, whereas the dashed boundary is exclusive. This example highlights the complexity involved in determining optimal policies.

A similar monotonicity result holds for the parameters of the compound Poisson process. We require the definition of the usual stochastic order that quantifies the concept of one random variable being ``smaller'' than another random variable.

\begin{definition}[The usual stochastic order]
Let $X$ and $Y$ be two random variables such that $\mathbb{P}(X > x) \leq \mathbb{P}(Y > x)$ for all $x \in \mathbb{R}$. Then $X$ is said to be smaller than $Y$ in the usual stochastic order, denoted by $X \leq_\textrm{st} Y$~\citep[1.A.1]{shaked2007stochastic}.
\end{definition}

Applying the concept of the usual stochastic order, we establish the theorem on monotonicity in degradation parameters.

\begin{theorem}[Monotonicity in degradation parameters]\label{thm:degparams}
Let $x = (x_1, \lambda_1, \phi_1, \ldots, x_M, \lambda_M, \phi_M)$ denote a state. Let $\mathcal{M}(x) \subseteq \mathcal{M}$ denote the set of machines for which it is optimal to do maintenance in state $x$. For any $x^\prime = (x_1, \lambda^\prime_1, \phi^\prime_1, \ldots, x_M, \lambda^\prime_M, \phi^\prime_M)$ that represents a state with worse degradation parameters than $x$, i.e.,
\begin{enumerate}
 \item $\lambda^\prime_m \geq \lambda_m$ and $\phi^\prime_m \geq \phi_m$ for all $m \in \mathcal{M}(x)$, and
 \item $\lambda^\prime_m = \lambda_m$ and $\phi^\prime_m = \phi_m$ otherwise.
\end{enumerate}

If $\phi^\prime_m \geq \phi_m$ implies that the corresponding shock size distributions satisfy $Y^\prime_m \geq_\textrm{st} Y_m$ (where $Y^\prime_m$ and $Y_m$ are the random variables associated with the respective shock size distributions), then it holds that
\begin{equation*}
 a^*(x) \in \textrm{arg} \underset{a\in\mathcal{U}(x^\prime)}{\textrm{min}} Q(x^\prime,a).
\end{equation*}
\end{theorem}

\mbox{Theorem \ref{thm:degparams}} essentially states that if a state $x^\prime$ is ``worse'' than another state $x$ in terms of degradation parameters, then the optimal maintenance action remains the same. Specifically, if the degradation parameters $\lambda^\prime_m$ in state $x^\prime$ are greater than or equal to $\lambda_m$ in state $x$, and the shock size distributions $Y^\prime_m$ are larger than $Y_m$ in the usual stochastic order, then the same maintenance actions should be applied in both states. Note that \mbox{Theorem \ref{thm:degparams}} applies to a broad range of MDPs with shock size distributions from the one-parameter exponential family, particularly the geometric distribution supported on $\mathbb{Z}_+$. 

While solving for an optimal policy explicitly can be computationally infeasible for larger asset networks, the derived structural properties provide valuable insights that can guide policy design. Although trained DRL agents generally do not satisfy structural guarantees such as the monotonicity properties derived here, the results remain valuable in broader DRL contexts, as they provide theoretical guidance for designing effective heuristic policies and for reducing the complexity of the DRL algorithm by shrinking the action space.
In DRL, such heuristic strategies play a crucial role: they serve not only as benchmarks for evaluating trained agents but also as initialization strategies that help kickstart learning in challenging environments. \Citet{kdtmpa}, for example, use carefully designed heuristic policies to initialize the API algorithm and demonstrate that such initialization significantly reduces the number of iterations needed to reach a near-optimal solution.

\section{Heuristic solution approaches}\label{sec:heuristics}

Heuristic methods offer practical and efficient strategies for solving complex decision-making problems, particularly when exact solutions are computationally infeasible. In this section, we discuss three heuristic approaches: the two-threshold control limit heuristic, the integrated Bayes heuristic and approximate policy iteration for BMDPs. We will now define each of these approaches.

\subsection{Two-threshold control limit heuristic}\label{subsec:control_limit}
The two-threshold control limit heuristic is an offline heuristic approach requiring at least information level $\infolevel{1}$. This heuristic involves two distinct thresholds: one for initiating preventive maintenance (PM) and another for opportunistic preventive maintenance (OPM). The PM threshold $\tau_m^\textrm{PM} \in (0, \xi_m]$ triggers scheduled maintenance activities when the $m$-th asset's degradation signal is greater than or equal to $\tau_m^\textrm{PM}$. Meanwhile, the OPM threshold $\tau_m^\textrm{OPM} \in (0, \tau_m^\textrm{PM}]$ is set at a lower degradation level, prompting maintenance activities when another asset is already undergoing maintenance. This approach capitalizes on the opportunity to address potential issues while sharing maintenance setup costs, thereby reducing the overall costs. 
We refer to this approach as $\pi^{\infolevel{1}}_\mathcal{N}$ because it is the most naive heuristic approach; it does not explicitly consider component heterogeneity, yet it is widely used in practice. Note that setting $\tau_m^\textrm{OPM} = \tau_m^\textrm{PM} \geq \xi_m$ for all $m\in \mathcal{M}$ yields a reactive heuristic, which we shall denote by $\pi^{\infolevel{1}}_\mathcal{R}$.

\subsection{Integrated Bayes heuristic}\label{subsec:integrated_bayes}
 
The integrated Bayes heuristic approach is adopted from~\citet[Section 5.3]{msom}, and applies a state-space truncation to solve the objective of \mbox{Eq. (\ref{eq:objectiveBMDP})} in the single asset scenario using value iteration. We extend the policy to the multi-asset scenario by applying the obtained policy $\pi^{\infolevel{1}}_{\mathcal{I}}$ to each asset individually, thus ignoring the economic dependence. Due to the curse of dimensionality, solving the objective of \mbox{Eq. (\ref{eq:objectiveBMDP})} for larger asset networks is computationally intractable. 

\subsection{Approximate policy iteration for Bayesian Markov decision processes}\label{subsec:api}

To learn maintenance policies for economically coupled assets in the context of BMDPs, we employ a variation of approximate policy iteration. Specifically, we utilize deep controlled learning (DCL), which was first introduced by~\citet{temizoz2023deep}. \Citet{kdtmpa} have shown that API/DCL can effectively learn maintenance policies (together with dispatching strategies) for an industrial-scale network of homogeneous assets, and that it outperforms several baseline heuristics, including decomposition-based methods and solutions derived from combinatorial optimization algorithms. Moreover, DCL has demonstrated superior performance in inventory problems compared to other DRL algorithms including \emph{proximal policy optimization} and \emph{asynchronous advantage actor-critic}~\citep{temizoz2023deep}.

To apply DCL to BMDPs, we introduce various novel techniques: (i) randomized action selection to avoid indexation bias, (ii) sampling from the prior distributions to directly train neural network policies for BMDPs, (iii) leveraging heuristic information to limit the action space in a significant number of states, and (iv) an open-loop feedback feature vector that enables the application of neural network policies trained in an $\infolevel{2}$ setting to be applied in an $\infolevel{1}$ setting.

\subsubsection{Deep controlled learning for Bayesian Markov decision processes}\label{subsubsec:DCL_details}
Following the approach of~\citet{kdtmpa}, we apply DCL to train a neural network for sequential action selection in asset management, thereby significantly reducing the complexity of the action space. To avoid indexation bias, we randomly select a new permutation $\sigma \in S(\mathcal{M})$, where $S(\mathcal{M})$ denotes the set of all permutations of $\mathcal{M}$, at each decision epoch, ensuring that the order in which actions are chosen for machines is randomized. Following each action selection, the input state is updated with the effects of the action before the algorithm decides on an action for the next asset. More specifically, the state transition $h \rightarrow h^\prime$ is decomposed into two stages: The first stage is governed by the deterministic outcomes of the selected actions $a \in \mathcal{U}(h)$, while the second stage is influenced by the random progression of the degradation processes. In more detail, $h \rightarrow h^\prime$ is decomposed into $h \overset{a_{\sigma{(1)}}}{\rightarrow} h^{a_{\sigma{(1)}}} \overset{a_{\sigma{(2)}}}{\rightarrow} \ldots \overset{a_{\sigma{(M)}}}{\rightarrow} h^{a_{\sigma{(M)}}} =: h^a$ and to $h^a \overset{t \rightarrow t+1}\longrightarrow h^\prime$. The sequence in which the actions $a_1, \ldots, a_M$ are handled does not matter, so we will only describe the processing of the action for the $m$-th asset.

In the case that $(a_m = 0)$, maintenance on the asset is postponed and the corresponding state variables remain unchanged. The action $(a_m= 1)$ represents initiating a maintenance action. The corresponding degradation level $x_m$ is set to $0$ and new degradation process parameters $\lambda_m$ and $\phi_m$ are sampled from their respective sampling distributions. To determine $h'$, we simply update the state variables of each asset based on the random evolution of the degradation processes.

We now provide a concise overview of the application of DCL to BMDPs (on the application of DCL to MDPs, we refer to~\citet{kdtmpa} and to~\citet{temizoz2023deep}). Under the information level $\infolevel{}$, the DCL algorithm comprises the following three steps:
\begin{enumerate}
\item Select an appropriate initial solution $\pi^{\infolevel{}}_0$.
\item Using $\pi^{\infolevel{}}_0$, construct a data set $\mathcal{D}$ containing state-action mappings.
\item Train a neural network classifier to learn the state-action mappings in $\mathcal{D}$.
\end{enumerate}

In step three, the neural network can be regarded as a parameterized function mapping from $\mathbb{R}^r$ to $\mathbb{R}^s$ for some $r,s \in \mathbb{N}$. We denote this function as $N_\theta(\cdot)$, where $\theta$ signifies the function parameters. The input to the neural network is the feature representation $f^{\infolevel{}}(h)\in \mathbb{R}^r$ of state $h$, and the output $N_\theta(\cdot) \in \mathbb{R}^{s}$, where $s=2$ in our case, is converted into a probability distribution over the action space. The action $\tilde{a}$ with the highest probability, $N_\theta( \cdot )_{\tilde{a}}$, is selected, effectively defining the policy. Given the actions $a_{\sigma(1)}, \ldots, a_{\sigma(m-1)}$ for the first $m-1$ assets (according to the permutation $\sigma$), the action $a_{\sigma(m)}$ for the asset $\sigma(m)$ in state $h$ is selected using the following decision rule:
\begin{equation*}
 \pi^{\sigma(m)}_{\theta}( f^{\infolevel{}}(h^{ a_{\sigma(m-1)}}) ) = \underset{\tilde{a}\in \mathcal{U}_m(h^{ a_{\sigma(m-1)} })}{\textrm{arg max}} [(N_\theta( f^{\infolevel{}}(h^{ a_{\sigma(m-1)} })) )_{\tilde{a}}],
\end{equation*}
where by convention $h^{ a_{\sigma(0)}} = h$. We refer to $\pi^{\infolevel{}}_{\theta}$ as the neural network policy, trained in a setting with information level $\infolevel{}$, that selects in each decision epoch, for each asset $\sigma(m)$, the action $\pi^{\sigma(m)}_{\theta}(f^{\infolevel{}}(h^{ \pi^{ \sigma(m-1)}_{\theta}}))$.

To apply DCL to BMDPs, we construct a data set $\mathcal{D}$ containing state-action mappings for BMDP states $\tilde{h}$. To advance these states, we modify the two-stage decomposition as follows: The first stage, $\tilde{h} \rightarrow \tilde{h}^a$, remains largely unchanged, except that after a maintenance action $(a_m = 1)$ is selected, the prior distributions $\tilde{\Lambda}_m$ and $\tilde{\Phi}_m$ are reset by setting the state variables $k_m(t)$ and $t_m(t)$ to 0. To determine $\tilde{h}^a \overset{t \rightarrow t+1}\longrightarrow \tilde{h}^\prime$, we first sample the degradation process parameters from the prior distributions specified in \mbox{Eq. (\ref{eq:prior_distr})}. Using these parameters, we sample the shock arrivals and the resulting increase in degradation, and subsequently update the state variables and belief distributions through the prior-to-posterior update. The posterior distributions are then used as the prior distributions for the next time step. This enables us to directly train policies using DCL in the BMDP setting.

In the next section, we discuss suitable initial policies and feature representations for state information.

\subsubsection{Initial policies and feature representations}\label{subsubsec:feature_representation}

\citet[Section 5.3]{kdtmpa} argue that initiating DCL with an appropriate policy $\pi^{\infolevel{}}_0$ considerably reduces computation times. The optimized two-threshold control limit heuristic satisfies all relevant listed properties since it explores sufficiently many states and has low computational complexity (as opposed to the integrated Bayes heuristic from \mbox{Section \ref{subsec:integrated_bayes}}). Moreover, the optimized PM and OPM thresholds $\tau_m^\textrm{PM}$ and $\tau_m^\textrm{OPM}$ can be leveraged to effectively restrict the individual, state-dependent action sets $\mathcal{U}_m(h)$ for a significant number of states. This is particularly useful when the failure threshold $\xi_m$ is relatively large. Specifically, we train and evaluate the policy improvements of the two-threshold control limit on a restricted action space defined as $\tilde{\mathcal{U}}(h) = \tilde{\mathcal{U}}_1(h) \times \ldots \times \tilde{\mathcal{U}}_M(h)$, where for each $m \in \mathcal{M}$, the individual action set $\tilde{\mathcal{U}}_m(h)$ is given by:
\begin{equation*}
 \tilde{\mathcal{U}}_m(h) = 
 \begin{cases}
 \{0\} & \text{if } x_m \leq \delta_m\cdot\tau_m^\textrm{OPM}, \\
 \{1\} & \text{if } x_m \geq \zeta_m\cdot\tau_m^\textrm{PM}, \\
 \{0,1\} & \text{otherwise.}
 \end{cases}
\end{equation*}
Here, for each $m \in \mathcal{M}$, $\delta_m \in [0,1]$ and $\zeta_m \geq 1$ are chosen conservatively to ensure that the action space is restricted only in clearly suboptimal regions of the state space. 
In all numerical experiments, we set $\delta_m \equiv 0.5$ and $\zeta_m \equiv 1.5$, so that PM is prohibited only when $x_m$ is well below the OPM threshold, and doing nothing is prohibited only when $x_m$ is substantially above the PM threshold.

The feature representation of a state depends on the information available to the manager. We propose a feature design that communicates the state information for each asset in a manner that is \emph{customized to the specific asset for which we are currently determining an action}. \Citet[Section 7.1]{kdtmpa} demonstrate that using such a feature design results in less training variability and consistently produces policies that perform significantly better.

In case of the underlying MDP under information level $\infolevel{2}$, the most compact representation of the state $h$ is given by
\begin{align*}
f^{\infolevel{2}}_1(h) = \left( x_1, \lambda_1, p_1, \iota_1, \ldots, x_M, \lambda_M, p_M, \iota_M, \eta \right),
\end{align*}
that is, the feature vector $f^{\infolevel{2}}_1(h)$ includes an information block $(x_m, \lambda_m, p_m, \iota_m)$ for each $m\in\mathcal{M}$ that can be derived from $h$, along with one additional feature. Here, $x_m$ is the observed degradation level of the $m$-th asset, and the entries $\lambda_m$ and $p_m$ denote its degradation parameters (cf. Example~\ref{example:geo}). The last block entry $\iota_m$ indicates whether we are currently selecting an action for asset $m$. Lastly, the additional feature $\eta$ indicates whether there is an opportunity for OPM. More specifically, $\eta$ equals $1$ if maintenance on at least one asset is mandatory or if the maintenance action has already been selected in the current decision epoch, and $0$ otherwise. Note that the dimension of the feature vector is $r = 4M+1$. 

In the case that we are agnostic about the degradation parameters, these features need to be estimated from the available history. This can be done via maximum likelihood estimation or by collapsing the belief distribution into a point estimate. Note that the maximum likelihood estimators can only be computed when some data has been collected. Bayesian inference does not suffer from this drawback since the knowledge of the true hyperparameters already yields a good estimate. The belief distribution can be collapsed in a suitable feature representation using an open-loop feedback approach. Specifically, we collapse the belief distributions into the following point estimates, which represent the means of the distributions:
\begin{equation}\label{eq:bayes_estimates}
 (\lambda^\textrm{Bayes}_m, p^\textrm{Bayes}_m) = \left( \frac{\alpha_m(t)}{ \beta_m(t) }, \frac{r_m(t)}{r_m(t) + s_m(t)} \right),
\end{equation}
where $\alpha_m(t), \beta_m(t), r_m(t)$ and $s_m(t)$ are defined in \mbox{Eqs. (\ref{eq:alpha_update})--(\ref{eq:b_update})}.

To apply an $\infolevel{2}$ trained policy in the $\infolevel{1}$ setting, we modify the feature representation $f^{\infolevel{2}}_1(h)$ by substituting the estimates from \mbox{Eq. (\ref{eq:bayes_estimates})}, i.e.,
\begin{align*}
f^{\infolevel{1}}_2(\tilde{h}) = \left( x_1(t), \lambda^\textrm{Bayes}_1, p^\textrm{Bayes}_1, \iota_1, \ldots, x_M(t), \lambda^\textrm{Bayes}_M, p^\textrm{Bayes}_M, \iota_M, \eta \right).
\end{align*}
We denote the resulting policy as $\pi^{\infolevel{2}}_{\theta}(f_2^{\infolevel{1}}(\tilde{h}))$ to emphasize that the policy is trained in the $\infolevel{2}$ setting but applied in the $\infolevel{1}$ setting using the open-loop feedback approach.

Lastly, we extend this feature representation for BMDPs by including the relevant history captured by $k_m(t)$ (the number of shocks) and $t_m(t)$ (the machine age). Thus, 
\begin{equation*} 
f^{\infolevel{1}}_3(\tilde{h}) = \left( x_1(t), \lambda^\textrm{Bayes}_1, p^\textrm{Bayes}_1, k_1(t), t_1(t), \iota_1, \ldots, x_M(t), \lambda^\textrm{Bayes}_M, p^\textrm{Bayes}_M, k_M(t), t_M(t), \iota_M, \eta \right).
\end{equation*}

It is important to note that these feature designs can be effective even outside the Bayesian framework assumed so far and are not constrained by the modeling assumptions inherent to that framework. While their accuracy is generally highest when the data exhibits behavior similar to the underlying model structure, this similarity is not a strict requirement for their practical usefulness.
\section{Simulation study}\label{sec:sim_study}

This section presents the findings of a compact simulation study where we optimize the decision process for two example instances. Unlike the case study that will be presented in \mbox{Section \ref{sec:case_study}}, this simulation study starts from the premise that the manager has full distributional information about the underlying MDP model, i.e., information level $\infolevel{1}$ or above. This allows for a strictly controlled setting. 

The aim of this simulation study is twofold:
\begin{enumerate}
 \item To assess the benefits of integrating learning and decision-making, which explicitly takes into account the heterogeneity in asset degradation (value of integration).
 \item To examine the value of acquiring additional information about the system to improve decision-making, particularly in understanding the uncertainty and variability of asset degradation (value of information).
\end{enumerate}

We compare the performance of each proposed solution approach with that of a policy optimized under information level $\infolevel{2}$. We restrict our analysis to the case where $M=2$ as it facilitates visualization of the trained policies. Furthermore, we assume that all component replacements come from the same pool of components, meaning that $\Lambda_m \overset{\textrm{d}}{\equiv} \Lambda$ and $\Phi_m \overset{\textrm{d}}{\equiv} \Phi$. For each instance of the simulation study, the true hyperparameters of the gamma distribution $\Lambda$ and the beta distribution $\Phi$ that model the population heterogeneity are denoted by $\alpha$, $\beta$, $r$ and $s$, and are listed in \mbox{Table \ref{tab:instance_parameters}}. The means of the distributions $\Lambda$ and $\Phi$ are denoted by $\mu_{\Lambda}$ and $\mu_{\Phi}$, respectively. These instances are selected to reflect increasing volatility in component heterogeneity, represented by an increasing coefficient of variation (CV) of the distributions, as well as a range of representative cost parameters. 

\begin{table}[!ht]
\centering
\resizebox{\textwidth}{!}{$
\begin{tabular}{ccc|cccc|cccc|cccc}
\toprule
Instance & $M$ & $\xi_m$ & $\mu_{\Lambda}$ & $\textrm{CV}_{\Lambda}$ & $\mu_{\Phi}$ & $\textrm{CV}_{\Phi}$ & $\alpha$ & $\nicefrac{1}{\beta}$ & $r$ & $s$ & $c_m^{\textrm{PM}}$ & $c_m^{\textrm{CM}}$ & $c^{\textrm{ST}}$ & $\gamma$ \\ \midrule
\textrm{I.1} & $2$ & $20$ & $1$ & $0.3$ & $0.5$ & $0.01$ & $11.11$ & $0.09$ & $4999.5$ & $4999.5$ & 1 & 5 & 1 & 0.99 \\
\textrm{I.2} & $2$ & $20$ & $1$ & $0.6$ & $0.5$ & $0.02$ & $2.78$ & $0.36$ & $1249.5$ & $1249.5$ & 1 & 10 & 1 & 0.99\\
\bottomrule
\end{tabular}
$}
\caption{Hyperparameter settings and cost structures (adopted from~\citet[Section 5]{msom}) considered in the simulation study.}
\label{tab:instance_parameters}
\end{table}

All performance results in this section are retrieved using $10^6$ repetitions of length $10^3$ time units. The reported half-widths represent asymptotic $95\%$ confidence intervals and account solely for variability within the model.

\subsection{Performance analysis of heuristic solution approaches}\label{subsec:sim_study_heuristics}

\begin{figure}[!ht]
 \centering
 \begin{subfigure}[c]{0.47\textwidth}
 \includegraphics[width=\linewidth]{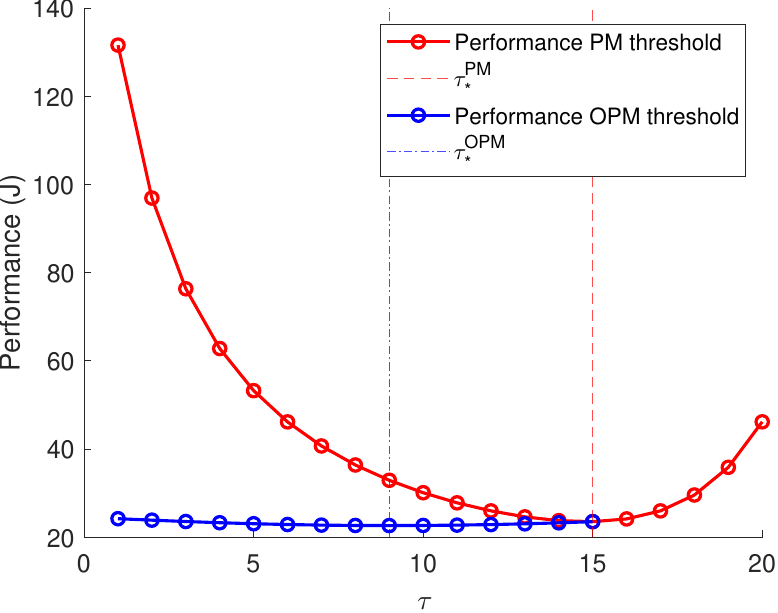} 
 \subcaption{The best-found PM and OPM thresholds for instance I.1 are $\tau_*^\textrm{PM} = 15$ and $\tau_*^\textrm{OPM} = 9$, respectively.}
 \label{subfig:threshold_optimization_test_I.1}
 \end{subfigure}\hfill
 \begin{subfigure}[c]{0.47\textwidth}
 \includegraphics[width=\linewidth]{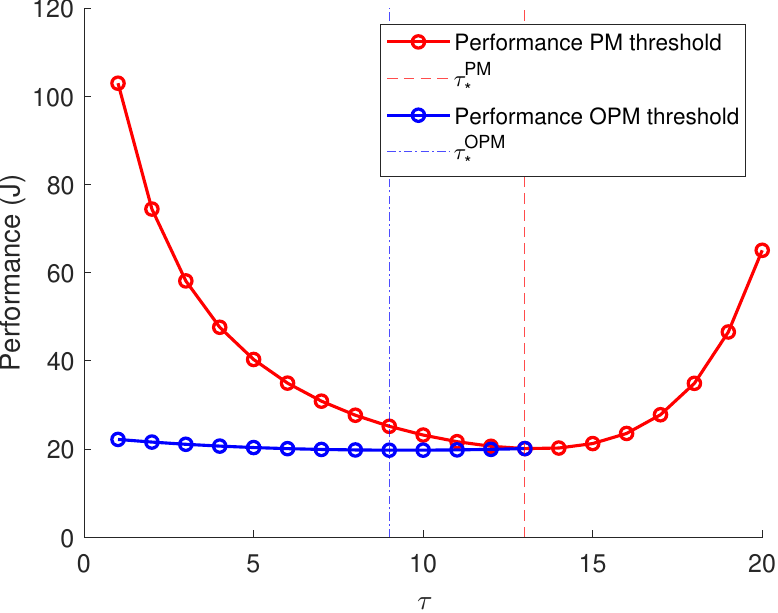} 
 \subcaption{The best-found PM and OPM thresholds for instance I.2 are $\tau_*^\textrm{PM} = 13$ and $\tau_*^\textrm{OPM} = 9$, respectively.}
\label{subfig:threshold_optimization_test_I.2}
 \end{subfigure}
 \caption{Results from the two-step simulation-based optimization of the two-threshold control limit heuristic for the instances I.1 and I.2.}
\label{fig:threshold_optimization_test_instances}
\end{figure}

For the instances of the simulation study, the two-threshold control limit (see \mbox{Section \ref{subsec:control_limit}}) is optimized using a two-step simulation-based process. Given that all component replacements draw from the same pool of components, we can exploit this symmetry by assuming identical PM and OPM thresholds for each asset, i.e., $\tau_m^\textrm{PM} \equiv \tau^\textrm{PM}$ and $\tau_m^\textrm{OPM} \equiv \tau^\textrm{OPM}$. Initially, the PM threshold is optimized without OPM. Once this optimal PM threshold $\tau_*^\textrm{PM}$ is established, we focus on determining the optimal OPM threshold $\tau_*^\textrm{OPM}$ under the PM threshold $\tau_*^\textrm{PM}$. This sequential approach not only simplifies the optimization process but also significantly reduces computational complexity. The optimization procedure for both instances I.1 and I.2 is visualized in \mbox{Figure \ref{fig:threshold_optimization_test_instances}}. As anticipated, increasing corrective maintenance costs combined with greater variability in degradation process parameters results in a lower PM threshold $\tau_*^\textrm{PM}$. Optimizing the OPM threshold further reduces costs by $3.78\%$ for instance I.1 and $1.80\%$ for instance I.2. The performance results for the two-threshold control limit $\pi^{\infolevel{1}}_\mathcal{N}$, the reactive heuristic $\pi^{\infolevel{1}}_\mathcal{R}$ and the integrated Bayes heuristic $\pi^{\infolevel{1}}_\mathcal{I}$ are summarized in \mbox{Table \ref{tab:test_instance_control_limits}}. Notably, the two-threshold control limit outperforms the integrated Bayes heuristic, which performs poorly in the multi-asset scenario due to its disregard of economic dependence.

\begin{table}[!ht]
\centering
\begin{tabular}{c|ccccc}
\toprule
Instance & $\tau_*^\textrm{PM}$ & $\tau_*^\textrm{OPM}$ & $J(\pi^{\infolevel{1}}_\mathcal{N})$ & $J(\pi^{\infolevel{1}}_\mathcal{R})$ & $J(\pi^{\infolevel{1}}_\mathcal{I})$\\ \midrule
\textrm{I.1} & $15$ & $9$ & $22.645 \pm 0.007$ & $46.177 \pm 0.012$ & $24.715 \pm 0.007$ \\
\textrm{I.2} & $13$ & $9$ & $19.741 \pm 0.011$ & $65.069 \pm 0.031$ & $20.380 \pm 0.010$\\
\bottomrule
\end{tabular}
\caption{Summary of the performance results of the heuristic solution approaches for the instances I.1 and I.2.}
\label{tab:test_instance_control_limits}
\end{table}

\subsection{Improving heuristic solutions through approximate policy iteration}\label{subsec:test_instance_results}

For both instances I.1 and I.2, we improve the reactive heuristic $\pi^{\infolevel{1}}_\mathcal{R}$ and the two-threshold control limit heuristic $\pi^{\infolevel{1}}_\mathcal{N}$ by applying three policy improvement steps using DCL in both $\infolevel{1}$ and $\infolevel{2}$ settings. The training results for policies trained in the $\infolevel{2}$ setting (i.e., the underlying MDP) and their corresponding performance when applied in the $\infolevel{1}$ setting using the open-loop feedback approach are presented in \mbox{Table \ref{tab:result-L2}}. The training results for policies directly trained in the $\infolevel{1}$ setting (i.e., the BMDP) are presented in \mbox{Table \ref{tab:result-L1}}.

\begin{table}[!ht]
\centering
\resizebox{1.0\linewidth}{!}{%
\begin{tabular}{c|c|c|c|c|c}
\toprule
 & & \multicolumn{2}{c|}{I.1} & \multicolumn{2}{c}{I.2} \\ 
\hline
 \multirow{2}{*}{\rotatebox[origin=c]{90}{\footnotesize\centering Gen 0}} & $\pi_0^{\infolevel{}}$ & $\pi^{\infolevel{1}}_\mathcal{N}$ & $\pi^{\infolevel{1}}_\mathcal{R}$ & $\pi^{\infolevel{1}}_\mathcal{N}$ & $\pi^{\infolevel{1}}_\mathcal{R}$ \\ 
\cline{2-6}
 & $J(\pi_0^{\infolevel{} })$ & $22.645 \pm 0.007$ & $46.177 \pm 0.012$ & $19.741 \pm 0.011$ & $65.069 \pm 0.031$ \\ \cline{1-6}
\multirow{2}{*}{\rotatebox[origin=c]{90}{\footnotesize\centering Gen 1}} & \cellcolor{gray!50}$J\big( \pi^{\infolevel{2}}_{\theta_{1}} (f^{\infolevel{2}}_1(h)) \big)$ & \cellcolor{gray!50}$21.585 \pm 0.007$ & \cellcolor{gray!50}$34.225 \pm 0.010$ & \cellcolor{gray!50}$18.487 \pm 0.010$ & \cellcolor{gray!50}$28.208 \pm 0.016$ \\ 
 & $J\big(\pi^{\infolevel{2}}_{\theta_{1}} (f^{\infolevel{1}}_2(\tilde{h}))\big)$ & $21.730 \pm 0.007$ & $36.522 \pm 0.010$ & $18.708 \pm 0.010$ & $ 38.188 \pm 0.023$ \\ \cline{1-6}
\multirow{2}{*}{\rotatebox[origin=c]{90}{\footnotesize\centering Gen 2}} & \cellcolor{gray!50}$J\big(\pi^{\infolevel{2}}_{\theta_{2}} (f^{\infolevel{2}}_1(h))\big)$ & \cellcolor{gray!50}$21.539 \pm 0.006$ & \cellcolor{gray!50}$25.138 \pm 0.008$ & \cellcolor{gray!50}$18.148 \pm 0.010$ & \cellcolor{gray!50}$20.523 \pm 0.012$ \\
 & $J\big(\pi^{\infolevel{2}}_{\theta_{2}} (f^{\infolevel{1}}_2(\tilde{h}))\big)$ & $21.725 \pm 0.007$ & $25.659 \pm 0.008$ & $18.435 \pm 0.010$ & $ 22.892 \pm 0.015$ \\ \cline{1-6}
\multirow{2}{*}{\rotatebox[origin=c]{90}{\footnotesize\centering Gen 3}} & \cellcolor{gray!50}$J\big(\pi^{\infolevel{2}}_{\theta_{3}} (f^{\infolevel{2}}_1(h))\big)$ & \cellcolor{gray!50}$\mathbf{21.518 \pm 0.006}$ & \cellcolor{gray!50}$24.258 \pm 0.007$ & \cellcolor{gray!50}$\mathbf{18.097 \pm 0.010}$ & \cellcolor{gray!50}$19.674 \pm 0.011$ \\
 & $J\big(\pi^{\infolevel{2}}_{\theta_{3}} (f^{\infolevel{1}}_2(\tilde{h}))\big)$ & $\mathbf{21.709 \pm 0.007}$ & $24.761 \pm 0.007$ & $\mathbf{18.376 \pm 0.010}$ & $21.514 \pm 0.013$ \\ 
\bottomrule
\end{tabular}
}
\caption{One-step policy improvement results for instances I.1 and I.2. \emph{\textcolor{darkgray}{Gray rows:}} The performance of the neural network policy $\pi^{\infolevel{2}}_{\theta_{}}$ in the $\infolevel{2}$ setting, trained on the underlying MDP. \emph{White rows:} The performance of the neural network policy $\pi^{\infolevel{2}}_{\theta_{}}$ applied in the $\infolevel{1}$ setting using the open-loop feedback approach. \textbf{Bold:} Indicates the lowest cost for each instance and information level across neural network generations.}
\label{tab:result-L2}
\end{table}

\begin{table}[!ht]
\centering
\resizebox{1.0\linewidth}{!}{%
\begin{tabular}{c|c|c|c|c|c}
\toprule
 & & \multicolumn{2}{c|}{I.1} & \multicolumn{2}{c}{I.2} \\ 
\hline
\multirow{2}{*}{\rotatebox[origin=c]{90}{\tiny\centering Gen 0}} & $\pi_0^{\infolevel{}}$ & $\pi^{\infolevel{1}}_\mathcal{N}$ & $\pi^{\infolevel{1}}_\mathcal{R}$ & $\pi^{\infolevel{1}}_\mathcal{N}$ & $\pi^{\infolevel{1}}_\mathcal{R}$ \\ 
\cline{2-6}
 & $J(\pi_0^{\infolevel{}})$ & $22.645 \pm 0.007$ & $46.177 \pm 0.012$ & $19.741 \pm 0.011$ & $65.069 \pm 0.031$ \\ \hline
\rotatebox[origin=c]{90}{\tiny\centering Gen 1} & $J\big(\pi^{\infolevel{1}}_{\theta_{1}}(f^{\infolevel{1}}_3(\tilde{h}))\big)$ & $21.834 \pm 0.007$ & $35.079 \pm 0.010$ & $18.611 \pm 0.010$ & $30.948 \pm 0.018$ \\ \hline
\rotatebox[origin=c]{90}{\tiny\centering Gen 2} & $J\big(\pi^{\infolevel{1}}_{\theta_{2}}(f^{\infolevel{1}}_3(\tilde{h}))\big)$ & $21.729 \pm 0.007$ & $25.795 \pm 0.008$ & $\mathbf{18.335 \pm 0.010}$ & $21.157 \pm 0.012$ \\ \hline
\rotatebox[origin=c]{90}{\tiny\centering Gen 3} & $J\big(\pi^{\infolevel{1}}_{\theta_{3}}(f^{\infolevel{1}}_3(\tilde{h}))\big)$ & $\mathbf{21.708 \pm 0.007}$ & $24.555 \pm 0.007$ & $18.340 \pm 0.010$ & $20.068 \pm 0.011$ \\
\bottomrule
\end{tabular}
}
\caption{One-step policy improvement results for instances I.1 and I.2, where the neural network policies $\pi^{\infolevel{1}}_{\theta_{}}(f^{\infolevel{1}}_3(\tilde{h}))$ are trained directly on the BMDP. \textbf{Bold:} Indicates the lowest cost for each instance across neural network generations.}
\label{tab:result-L1}
\end{table}

In conclusion, the results in \mbox{Tables \ref{tab:result-L2} and \ref{tab:result-L1}} demonstrate that initializing DCL with the two-threshold control limit $\pi^{\infolevel{1}}_\mathcal{N}$ significantly reduces the number of iterations needed to achieve near-optimal performance in both the information settings $\infolevel{1}$ and $\infolevel{2}$. The best-found neural network policies using the open-loop feedback approach $\pi^{\infolevel{2}}_{\theta_{}}(f^{\infolevel{1}}_2(\tilde{h}))$ achieve similar performance to the best-found neural network policies $\pi^{\infolevel{1}}_{\theta_{}}(f^{\infolevel{1}}_3(\tilde{h}))$ that are directly trained on the BMDP. The improvement of the best-found policy in setting $\infolevel{1}$ over the optimized two-threshold control limit is $4.14\%$ (I.1) and $7.12\%$ (I.2). This suggests that the value of integration is significant and increases with the volatility of component heterogeneity and the cost ratio $\nicefrac{c_m^\textrm{CM}}{c_m^\textrm{PM}}$. However, the performance difference between the best-found policy in settings $\infolevel{1}$ and $\infolevel{2}$ is only $0.88\%$ (I.1) and $1.32\%$ (I.2), suggesting that the value of information is relatively small when the parameter uncertainty is managed effectively. Finally, all trained policies outperform the integrated Bayes heuristic $\pi^{\infolevel{1}}_\mathcal{I}$, which is the current state-of-the-art, with the best-found policy improving on it by $12.17\%$ (I.1) and $10.03\%$ (I.2).

\begin{figure}[!ht]
 \centering
 \begin{subfigure}[t]{0.48\textwidth} \includegraphics[width=\linewidth]{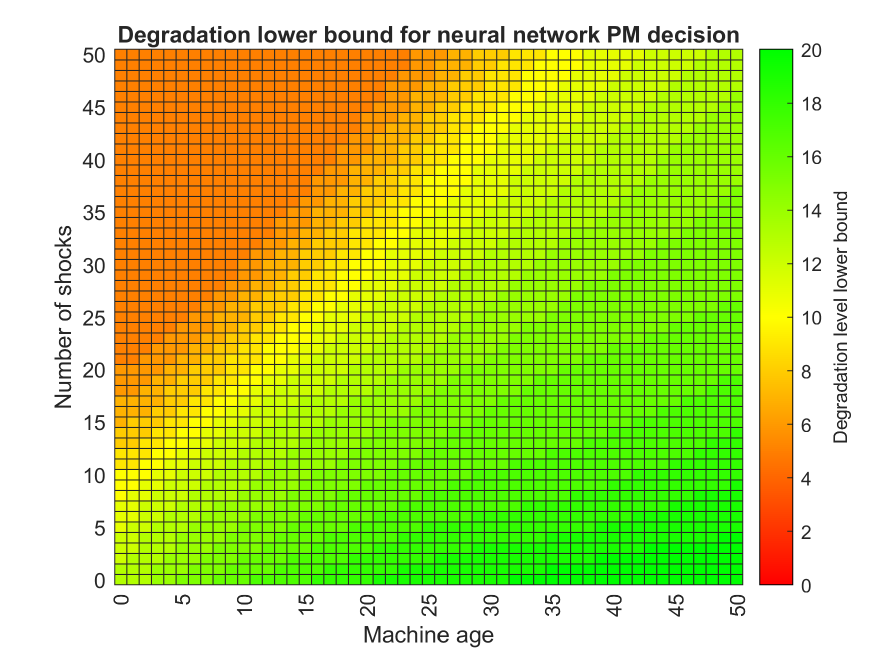} 
 \subcaption{Heatmap of the minimum degradation level $x_1$ at which maintenance is initiated, given $k_1(t)$ and $t_1(t)$, assuming the other machine is in the healthy state $\Big(x_2(t), k_2(t), t_2(t)\Big) = (0, 0, 0)$.}
 \label{subfig:neural_network_PM_decision_I.1}
 \end{subfigure}\hfill
 \begin{subfigure}[t]{0.48\textwidth}
 \includegraphics[width=\linewidth]{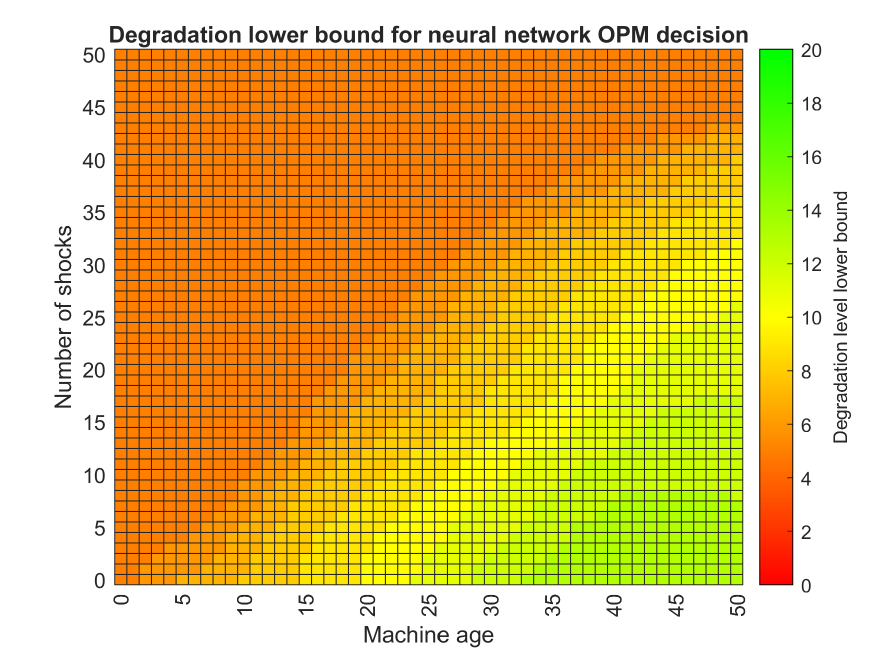} 
 \subcaption{Heatmap of the minimum degradation level $x_1$ at which maintenance is initiated, given $k_1(t)$ and $t_1(t)$, assuming the other machine is in the failed state $\Big(x_2(t), k_2(t), t_2(t)\Big) = \Big(\xi_2, \frac{\xi_2\cdot\mu_{\Phi}}{1-\mu_{\Phi}}, \frac{\xi_2\cdot\mu_{\Phi}}{(1-\mu_{\Phi})\cdot\mu_{\Lambda}}\Big)$.}
\label{subfig:neural_network_OPM_decision_I.1}
 \end{subfigure}
 \caption{Two policy slices from the best-performing neural network policy for instance I.1, illustrating the complex transformation from PM decisions to OPM decisions.}
\label{fig:neural_network_PM_OPM_decisions_I.1}
\end{figure}

Lastly, to illustrate the effectiveness of the DRL approach in learning OPM strategies, we present two policy slices in \mbox{Figure \ref{fig:neural_network_PM_OPM_decisions_I.1}} from the best-performing $\infolevel{1}$ neural network policy for instance I.1 (specifically, the 3rd-generation policy $\pi^{\infolevel{1}}_{\theta_{3}}(f^{\infolevel{1}}_3(\tilde{h}))$ improving the two-threshold control limit policy $\pi^{\infolevel{1}}_\mathcal{N}$). These slices demonstrate that DCL learns, in just a few iterations, a complex transformation from PM to OPM decisions\textemdash one that cannot be captured by a simple set of rules. Notably, from a managerial perspective, the minimum degradation level required for maintenance intervention generally decreases when there is an opportunity to share maintenance setup costs.

\section{Case study on degradation data of interventional X-ray system filaments}\label{sec:case_study}

\begin{figure}[!ht]
 \centering
 \begin{subfigure}[b]{0.5\textwidth}
\includegraphics[width=0.87\textwidth]{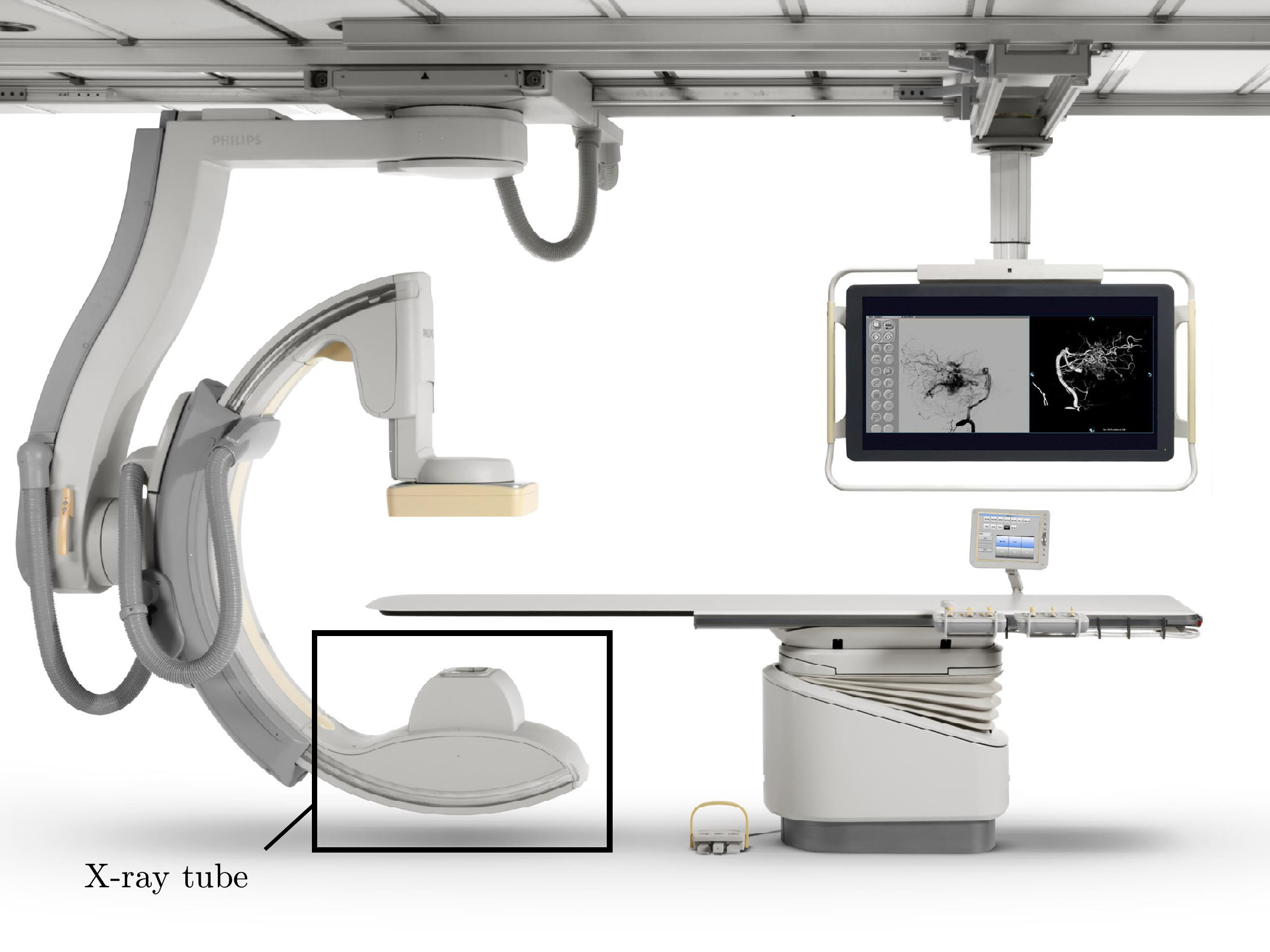}
\caption{An IXR system with the X-ray tube denoted by a rectangle.}
\label{ixrExample}
 \end{subfigure}\hfill
 \begin{subfigure}[b]{0.5\textwidth}
 \includegraphics[width=1\textwidth]{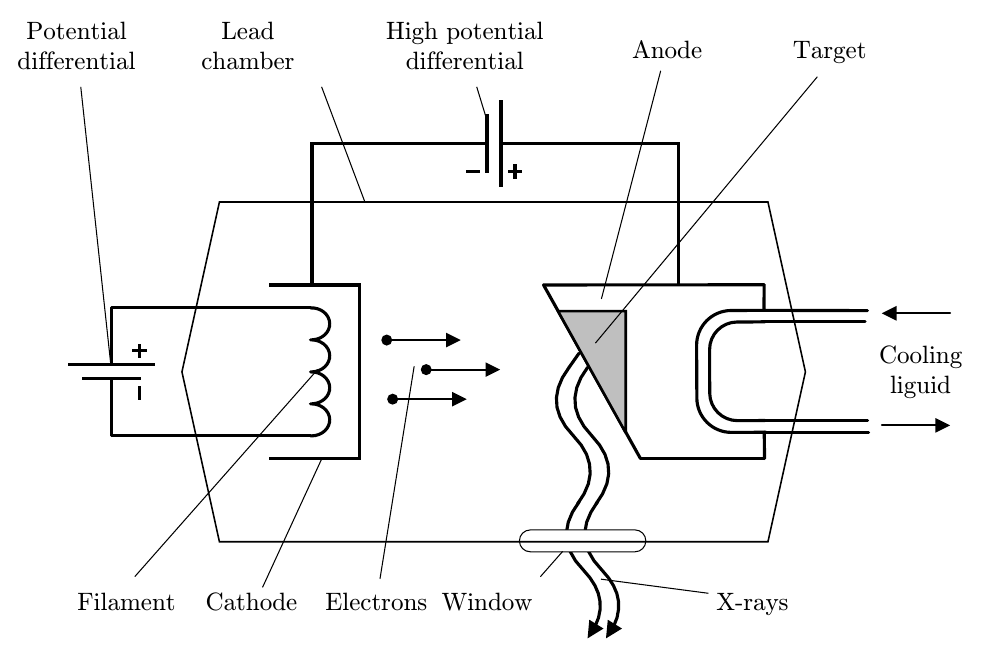}
\caption{Simplified X-ray tube schematic.}
\label{schematicXraytube}

 \end{subfigure}
 \caption{Diagram showing the positioning of the tungsten filament inside an IXR system.}
\label{fig:xray_filament}
\end{figure}

In this section, we assess the performance of the approaches outlined in \mbox{Section \ref{sec:heuristics}} on real-world degradation data of a component crucial to the functioning of medical imaging systems. Medical imaging devices, such as IXR systems, cost about one million USD, with annual maintenance expenses around $10\%$ of the initial cost~\citep{ECRI2013}. Over a typical 10-year lifespan, maintenance accounts for nearly half of the total ownership cost. X-ray tubes, the most expensive components of IXR systems, are critical for image-guided procedures but prone to failure, predominantly due to filament wear. During each imaging procedure, a high electric current is applied to the tungsten filament to generate electrons for the X-ray beam. This repeated heating causes small amounts of tungsten to evaporate, thinning the filament and forming a hotspot that eventually leads to failure~\citep{covington1973hot}. Importantly, degradation increments occur only during X-ray exposures and only when the electric current applied to the filament is sufficiently high to cause material degradation. Consequently, degradation accumulates through a sequence of operational events associated with system usage. In our modeling framework, these usage events are interpreted as the stochastic shocks that produce incremental degradation of the filament. We refer the reader to \mbox{Figure \ref{fig:xray_filament}} for a schematic representation of an IXR system.

To prevent unnecessary unavailability of such scanning equipment, Philips, which is a manufacturer of the IXR system, has developed a health indicator to monitor filament wear, collecting real-time degradation data from each IXR system. This generates a large data set containing real-time degradation data, recorded as time-series of filament wear, from the start of each filament's lifespan until either failure or the present, for all IXR systems in operation. The data set used is obtained from~\citet[Section 7]{msom} and originates from Philips. In addition to manufacturing, Philips offers maintenance and service contracts to hospitals using their IXR systems. Dutch hospitals typically operate several IXR scanners across various departments.

This case study serves as a realistic industrial illustration of the proposed methodology. The modeling framework and DRL-based solution approach are not specific to medical imaging systems; rather, they apply to asset networks in which degradation can be represented by stochastic increment or shock processes with learnable parameters, such as wind turbines, semiconductor equipment, or other high-value capital goods.

\subsection{Degradation data of X-ray tubes in IXR systems}\label{subsec:ixr_data}

The X-ray tube degradation data set consists of $52$ time series, each representing the degradation level of a distinct X-ray tube. Let $\mathcal{I}$ denote the set of X-ray tubes for which data is available, i.e., $|\mathcal{I}| = 52$. The time series for each X-ray tube $i \in \mathcal{I}$ is denoted as $ \mathcal{J}_i$. Each data point in the time series $\mathcal{J}_i$ is represented as a tuple $(t_j, x_j)_i$, where $t_j$ is the age of the X-ray tube, and $x_j$ is the degradation level at that age. Each tuple $(t_j, x_j)_i$ is generated when an IXR system is operated, and the time series contain between $20,000$ and $300,000$ data points, covering a time span of two to five years. These time series can be transformed into a series of BMDP trajectories, enabling us to directly evaluate the performance of our methods on the data.

For confidentiality reasons, the data was left-truncated and normalized. All time series start at $x_0 = 0$ for $t_0 = 0$ and end at $x_{|\mathcal{J}_i|} = 50$ (i.e., $\xi_m \equiv 50$ for all $i \in \mathcal{I}$). For each time series, the interarrival times \( (t_j - t_{j-1}) \) between consecutive data points and the degradation increments $(x_j - x_{j-1})$ were computed. To account for non-operational periods such as weekends, nights, or other extended downtimes, outliers in the interarrival times were removed. This allowed for the transformation of the original time series into those based on the operational age of each X-ray tube. Furthermore, several data points in the data where the image-guided procedure was deemed too short to cause significant wear on the X-ray tube were removed. Lastly, the time was normalized so that one unit of time corresponds roughly to the minimum operational time required for practical maintenance tasks, such as dispatching a service engineer to a hospital.

A statistical analysis of the resulting data set revealed that (i) there was no evidence to reject the assumption that shocks arrive according to a Poisson process, (ii) damage sizes are best modeled by a geometric distribution, and (iii) the data exhibits heterogeneity, meaning that the distribution parameters for interarrival times and shock sizes vary across components. Indeed, the estimated hyperparameters (see \mbox{Table \ref{tab:casestudy_parameters}}) indicate significant component heterogeneity. For the shock rate $\lambda_m$, the statistical analysis yields a mean of $1.414$ and a CV of $0.157$. Similarly, for the shock size parameter, the mean is $0.487$ with a CV of $0.234$.

\subsection{Numerical experiments and model calibration}\label{subsec:model_calibration}

In this section, we outline the model calibration process, that is, estimating the hyperparameters of the degradation model. The data set introduced in the previous section is representative of the typical conditions encountered by the asset manager in practice, i.e., the baseline information level $\infolevel{0}$. To address the heterogeneity of the component population, we estimate the hyperparameters of the distributions $\Lambda_m$ and $\Phi_m$ for each $m \in \mathcal{M}$ from the available X-ray tube degradation data using the maximum likelihood estimation procedure provided by~\citet[Online Appendix C]{msom}. As in \mbox{Section \ref{sec:sim_study}}, we assume that all component replacements stem from the same pool of components, i.e., $\Lambda_m \overset{\textrm{d}}{\equiv} \Lambda$ and $\Phi_m \overset{\textrm{d}}{\equiv} \Phi$. This assumption reflects asset-level homogeneity and is justified because, although the data may stem from distinct machines operated in different locations, the machines are of the same type, and the operating environment for each is controlled and standardized.

We divided the data set $\mathcal{I}$ into a training set $\mathcal{I_\text{train}}$ with $|\mathcal{I_\text{train}}| = 10$ and a test set $\mathcal{I_\text{test}}$ with $|\mathcal{I_\text{test}}| = 42$. The training set $\mathcal{I_\text{train}}$ is a randomly selected subset of $\mathcal{I}$ used solely for parameter estimation, while the test set $\mathcal{I_\text{test}}$ is reserved for evaluating our methods. For a detailed overview of the estimated parameters and cost settings used in this case study, refer to \mbox{Table \ref{tab:casestudy_parameters}}. The instance CS.1 features a single asset without consideration of economic dependence, which simplifies the problem to the maintenance problem studied in~\citet{msom}. In contrast, instances CS.2 and CS.3, which involve two and five machines respectively, more accurately reflect the complexity of a real-world hospital setting.

\begin{table}[!ht]
\centering
\begin{tabular}{cccc|cccc}
\toprule
 $\mu_{\Lambda}$ & $\textrm{CV}_{\Lambda}$ & $\mu_{\Phi}$ & $\textrm{CV}_{\Phi}$ & $\alpha$ & $\nicefrac{1}{\beta}$ & $r$ & $s$ \\ \midrule
$1.414$ & $0.157$ & $0.487$ & $0.234$ & $40.696$ & $28.779$ & $8.924$ & $9.405$ \\
\bottomrule
\end{tabular}
~\\
~\\
\centering
\begin{tabular}{ccc|cccc}
\toprule
Instance & $M$ & $\xi_m$ & $c_m^{\textrm{PM}}$ & $c_m^{\textrm{CM}}$ & $c^{\textrm{ST}}$ & $\gamma$ \\ \midrule
\textrm{CS.1} & $1$ & $50$ & 1 & 5 & 0 & 0.99 \\
\textrm{CS.2} & $2$ & $50$ & 1 & 5 & 1 & 0.99\\
\textrm{CS.3} & $5$ & $50$ & 1 & 5 & 1 & 0.99 \\
\bottomrule
\end{tabular}
\caption{Hyperparameter settings and cost structures considered in the case study.}
\label{tab:casestudy_parameters}
\end{table}

For the calibration of the heuristic solution approaches, we optimize both the two-threshold control limit heuristic and the integrated Bayes heuristic for the fitted degradation model. Note that the integrated Bayes approach is specifically optimized for the cost ratio of CS.1, where $\nicefrac{c_m^{\textrm{CM}}}{c_m^{\textrm{PM}}} = 5$, and the retrieved solution is used for the instances CS.2 and CS.3 as well. See \mbox{Table \ref{tab:case_study_instance_control_limits}} for a detailed summary of the heuristic solution calibration results. All performance results in this section are obtained from $10^6$ repetitions, each lasting $10^3$ time units. The reported half-widths again represent asymptotic $95\%$ confidence intervals.

We improve the calibrated two-threshold control limit using 3 iterations of DCL, for which we present the results in \mbox{Table \ref{tab:DCL_result-CS.1-3_L2}} and \mbox{Table \ref{tab:DCL_result-CS.1-3_L1}}.

\begin{table}[!ht]
\centering
\begin{tabular}{c|ccccc}
\toprule
Instance & $\tau_*^\textrm{PM}$ & $\tau_*^\textrm{OPM}$ & $J(\pi^{\infolevel{1}}_\mathcal{N})$ & $J(\pi^{\infolevel{1}}_\mathcal{R})$ & $J(\pi^{\infolevel{1}}_\mathcal{I})$ \\ \midrule
\textrm{CS.1} & $40$ & $-$ & $3.146 \pm 0.002$ & $11.071 \pm 0.006$ & $2.974 \pm 0.002$ \\
\textrm{CS.2} & $41$ & $28$ & $11.381 \pm 0.005$ & $26.516 \pm 0.010$ & $11.558 \pm 0.005$ \\
\textrm{CS.3} & $41$ & $29$ & $26.405 \pm 0.007$ & $65.876 \pm 0.016$ & $28.270 \pm 0.007$ \\
\bottomrule
\end{tabular}
\caption{Summary of the heuristic solution calibration results for the case study instances CS.1--3.}
\label{tab:case_study_instance_control_limits}
\end{table}

\begin{table}[!ht]
\centering
\begin{tabular}{c|c|c|c|c}
\toprule
 & & CS.1 & CS.2 & CS.3 \\ 
\hline
\multirow{2}{*}{\rotatebox[origin=c]{90}{ \footnotesize\centering Gen 0}} & $\pi_0^{\infolevel{}}$ & $\pi^{\infolevel{1}}_\mathcal{N}$ & $\pi^{\infolevel{1}}_\mathcal{N}$ & $\pi^{\infolevel{1}}_\mathcal{N}$ \\ 
\cline{2-5}
 & $J(\pi_0^{\infolevel{}})$ & $3.146 \pm 0.002$ & $11.381 \pm 0.005$ & $26.405 \pm 0.007$ \\ \cline{1-5}
\multirow{2}{*}{\rotatebox[origin=c]{90}{ \footnotesize\centering Gen 1}} & \cellcolor{gray!50}$J\big( \pi^{\infolevel{2}}_{\theta_{1}} (f^{\infolevel{2}}_1(h)) \big)$ & \cellcolor{gray!50}$2.932 \pm 0.002$ & \cellcolor{gray!50}$10.564 \pm 0.004$ & \cellcolor{gray!50}$24.020 \pm 0.006$ \\ 
 & $J\big(\pi^{\infolevel{2}}_{\theta_{1}} (f^{\infolevel{1}}_2(\tilde{h}))\big)$ & $3.936 \pm 0.003$ & $12.663 \pm 0.006$ & $28.498 \pm 0.008$ \\ \cline{1-5}
\multirow{2}{*}{\rotatebox[origin=c]{90}{ \footnotesize\centering Gen 2}} & \cellcolor{gray!50}$J\big(\pi^{\infolevel{2}}_{\theta_{2}} (f^{\infolevel{2}}_1(h))\big)$ & \cellcolor{gray!50}$2.90 \pm 0.002$ & \cellcolor{gray!50}$10.471 \pm 0.004$ & \cellcolor{gray!50}$23.660 \pm 0.006$ \\
 & $J\big(\pi^{\infolevel{2}}_{\theta_{2}} (f^{\infolevel{1}}_2(\tilde{h}))\big)$ & $\mathbf{3.767 \pm 0.003}$ & $\mathbf{12.527 \pm 0.006}$ & $\mathbf{28.185 \pm 0.008}$ \\ \cline{1-5}
\multirow{2}{*}{\rotatebox[origin=c]{90}{ \footnotesize\centering Gen 3}} & \cellcolor{gray!50}$J\big(\pi^{\infolevel{2}}_{\theta_{3}} (f^{\infolevel{2}}_1(h))\big)$ & \cellcolor{gray!50}$2.901 \pm 0.002$ & \cellcolor{gray!50}$10.516 \pm 0.004$ & \cellcolor{gray!50}$23.507 \pm 0.006$\\
 & $J\big(\pi^{\infolevel{2}}_{\theta_{3}} (f^{\infolevel{1}}_2(\tilde{h}))\big)$ & $3.804 \pm 0.003$ & $13.009 \pm 0.007 $ & $29.336 \pm 0.009$ \\ 
\bottomrule
\end{tabular}

\caption{One-step policy improvement results for case study instances CS.1--3. \emph{\textcolor{darkgray}{Gray rows:}} The performance of the neural network policy $\pi^{\infolevel{2}}_{\theta_{}}$ in the $\infolevel{2}$ setting, trained on the underlying MDP. \emph{White rows:} The performance of the neural network policy $\pi^{\infolevel{2}}_{\theta_{}}$ applied in the $\infolevel{1}$ setting using the open-loop feedback approach. \textbf{Bold:} Indicates the lowest cost for each instance under $\infolevel{1}$ across neural network generations.}
\label{tab:DCL_result-CS.1-3_L2}
\end{table}

\begin{table}[!ht]
\centering
\begin{tabular}{c|c|c|c|c}
\toprule
 & & CS.1 & CS.2 & CS.3 \\ 
\hline
\multirow{2}{*}{\rotatebox[origin=c]{90}{ \tiny\centering Gen 0}} & $\pi_0^{\infolevel{}}$ & $\pi^{\infolevel{1}}_\mathcal{N}$ & $\pi^{\infolevel{1}}_\mathcal{N}$ & $\pi^{\infolevel{1}}_\mathcal{N}$ \\ 
\cline{2-5}
 & $J(\pi_0^{\infolevel{}})$ & $3.146 \pm 0.002$ & $11.381 \pm 0.005$ & $26.405 \pm 0.007$ \\ \cline{1-5}
\hline
\rotatebox[origin=c]{90}{\tiny\centering Gen 1} & $J\big(\pi^{\infolevel{1}}_{\theta_{1}}(f^{\infolevel{1}}_3(\tilde{h}))\big)$ & $2.975 \pm 0.002$ & $11.003 \pm 0.005$ & $24.765 \pm 0.006$ \\ 
\hline
\rotatebox[origin=c]{90}{\tiny\centering Gen 2} & $J\big(\pi^{\infolevel{1}}_{\theta_{2}}(f^{\infolevel{1}}_3(\tilde{h}))\big)$ & $\mathbf{2.957 \pm 0.002}$ & $10.703 \pm 0.004$ & $\mathbf{24.197 \pm 0.006}$ \\ 
\hline
\rotatebox[origin=c]{90}{\tiny\centering Gen 3} & $J\big(\pi^{\infolevel{1}}_{\theta_{3}}(f^{\infolevel{1}}_3(\tilde{h}))\big)$ & $2.959 \pm 0.002$ & $\mathbf{10.601 \pm 0.004}$ & $24.605 \pm 0.006$ \\ 
\bottomrule
\end{tabular}
\caption{One-step policy improvement results for case study instances CS.1--3, where the neural network policies $\pi^{\infolevel{1}}_{\theta_{}}(f^{\infolevel{1}}_3(\tilde{h}))$ are trained directly on the BMDP. \textbf{Bold:} Indicates the lowest cost for each instance across neural network generations.}
\label{tab:DCL_result-CS.1-3_L1}
\end{table}

The findings are mostly consistent with the results presented in \mbox{Section \ref{sec:heuristics}}. Notable is that the performance of the neural network policies $\pi^{\infolevel{2}}_{\theta_{}}(f^{\infolevel{1}}_2(\tilde{h}))$ using the open-loop feedback approach is significantly worse compared to that of the neural network policies $\pi^{\infolevel{1}}_{\theta_{}}(f^{\infolevel{1}}_3(\tilde{h}))$ directly trained on the BMDP. This decline in performance is likely attributable to a combination of several factors: (i) the substantial increase in component heterogeneity, as reflected by the higher CV of the distributions $\Lambda$ and $\Phi$; (ii) the limited number of observations available prior to component replacement; (iii) the reduction of the posterior distribution to a single point estimate, which entails a considerable loss of information; and (iv) the risk that unrepresentative initial data biases the inferred distribution parameters, potentially leading the neural network policy to trigger premature interventions.

\begin{figure}[!ht]
 \centering
 \begin{subfigure}[t]{0.48\textwidth}
 \includegraphics[width=\linewidth]{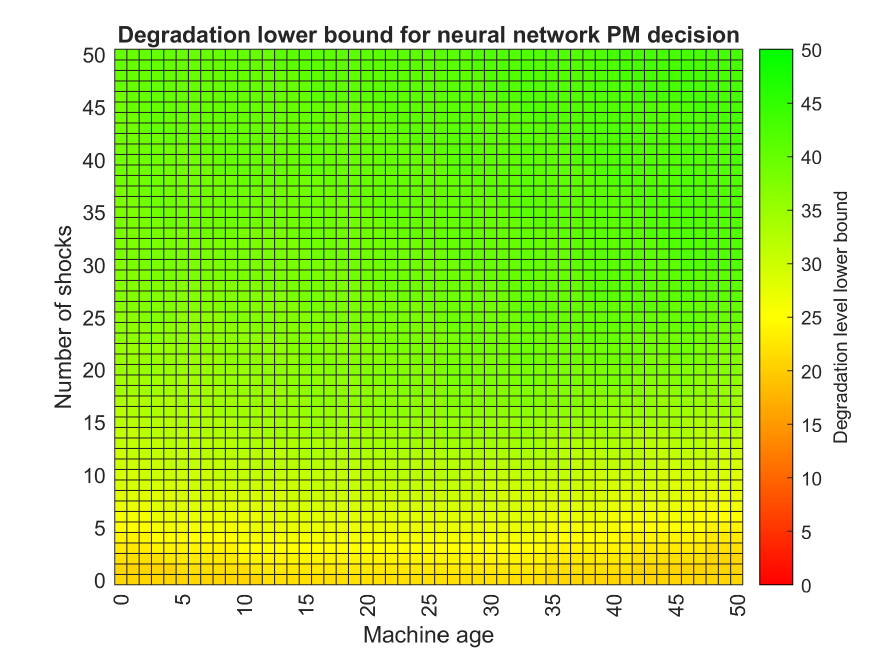} 
 \subcaption{Heatmap of the minimum degradation level $x_1$ at which maintenance is initiated according to the neural network policy, given $k_1(t)$ and $t_1(t)$.}
 \label{subfig:neural_network_PM_decision_CS.1}
 \end{subfigure}\hfill
 \begin{subfigure}[t]{0.48\textwidth}
 \includegraphics[width=\linewidth]{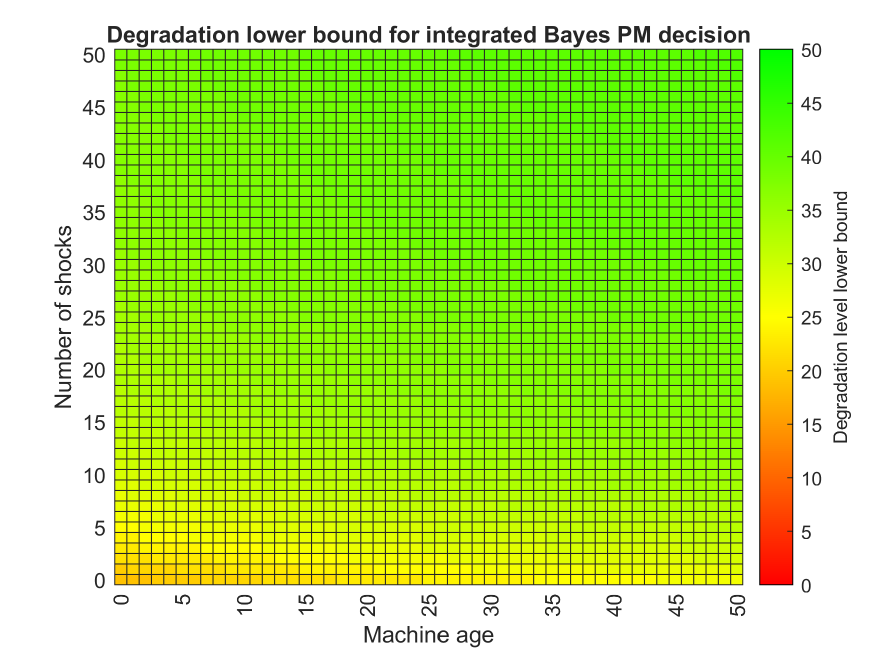} 
 \subcaption{Heatmap of the minimum degradation level $x_1$ at which maintenance is initiated according to the integrated Bayes heuristic, given $k_1(t)$ and $t_1(t)$.}
\label{subfig:IB_PM_decision_CS.1}
 \end{subfigure}
 \caption{Policy visualization of the best-performing neural network policy for instance CS.1 and the integrated Bayes heuristic, illustrating their similarity.}
\label{fig:neural_network_integrated_Bayes_PM_decisions_CS.1}
\end{figure}

For instance CS.1, the best-found neural network policy in the $\infolevel{1}$ setting achieves comparable performance to the integrated Bayes heuristic. Both policies are depicted in \mbox{Figure \ref{fig:neural_network_integrated_Bayes_PM_decisions_CS.1}} and exhibit a high degree of similarity.

An overlap analysis of \mbox{Figure \ref{fig:neural_network_integrated_Bayes_PM_decisions_CS.1}} demonstrates a strong alignment between the two policies. Overall, $13.23\%$ of the policy decisions were identical, and $65.36\%$ differed by no more than $2\%$ of the failure threshold $\xi_m$. These percentages increase significantly for the lower-left $20\times20$ subgrid, which contains more frequently visited states: $47.50\%$ of values were identical, and $97.25\%$ differed by at most $2\%$ of $\xi_m$. These results indicate that DCL can produce policies that closely approximate the near-optimal integrated Bayes heuristic for a real-world component, especially in regions of the state space that are frequently encountered. Moreover, DCL produces policies with state-of-the-art performance in a multi-asset setting that more closely resembles real-world hospital conditions.

\subsection{Numerical results and model validation}\label{subsec:model_validation}

In this section, we assess the predictive accuracy and robustness of the best solution approaches identified in \mbox{Tables \ref{tab:case_study_instance_control_limits}, \ref{tab:DCL_result-CS.1-3_L2} and \ref{tab:DCL_result-CS.1-3_L1}} by evaluating their performance directly on the test data. Specifically, we analyze the ability of the model to generalize beyond the training data and evaluate how well the estimated hyperparameters of the distributions $\Lambda$ and $\Phi$ capture the degradation behavior of other X-ray tubes across different case study instances. To this end, the time series in the test set $\mathcal{I_\text{test}}$ are transformed into BMDP trajectories. All performance results are derived from $10^4$ repetitions. The length of each repetition is set to $10^3$, with trajectories sampled from $\mathcal{I_\text{test}}$ with replacement. The reported half-widths correspond to asymptotic $95\%$ confidence intervals.

\begin{table}[!ht]
\centering
\resizebox{1.0\linewidth}{!}{%
\begin{tabular}{c|c|c|c|c|c}
\toprule
Instance & $J\big(\pi^{\infolevel{2}}_{\theta_{1}}(f^{\infolevel{1}}_2(\tilde{h}))\big)$ & $J\big(\pi^{\infolevel{1}}_{\theta_{1}}(f^{\infolevel{1}}_3(\tilde{h})) \big)$ & $J(\pi^{\infolevel{1}}_\mathcal{N})$ & $J(\pi^{\infolevel{1}}_\mathcal{R})$ & $J(\pi^{\infolevel{1}}_\mathcal{I})$ \\ \midrule
\textrm{CS.1} & $6.163 \pm 0.058$ & $4.278 \pm 0.046$ & $4.717 \pm 0.050$ & $11.270 \pm 0.080$ & $\mathbf{4.175 \pm 0.045}$ \\
\textrm{CS.2} & $16.623 \pm 0.099$ & $\mathbf{13.682 \pm 0.084}$ & $14.336 \pm 0.087$ & $26.956 \pm 0.135$ & $14.024 \pm 0.087$ \\
\textrm{CS.3} & $39.656 \pm 0.143$ & $\mathbf{31.444 \pm 0.123}$ & $33.507 \pm 0.128$ & $66.875 \pm 0.210$ & $34.288 \pm 0.135$ \\
\bottomrule
\end{tabular}
}
\caption{Model validation results for the case study instances CS.1--3. \textbf{Bold:} Highlights the lowest cost achieved on the test set for each case study instance across the proposed solution approaches.}
\label{tab:model_validation_case_study}
\end{table}

In summary, the model validation results presented in \mbox{Table \ref{tab:model_validation_case_study}} for case study instances CS.1--3 indicate that the fitted model generalizes well to unseen data. Integrating parameter estimation via the open-loop feedback approach results in a significant performance drop, as expected from the results in \mbox{Table \ref{tab:DCL_result-CS.1-3_L2}}. In contrast, the $\infolevel{1}$ DCL-improved policies consistently outperform the benchmark solutions across all case study instances, except for the single asset scenario, where they achieve near-optimal performance (relative to the fitted model parameters in \mbox{Table \ref{tab:casestudy_parameters}}).
\section{Conclusion}\label{sec:conclusion}
In this paper, we developed a novel maintenance model aimed at optimizing the management of a network of advanced industrial assets prone to costly and disruptive unplanned downtimes. By addressing the limitations of traditional condition-based maintenance models, which often assume homogeneous assets in isolation, our approach introduces a more practical framework that accounts for both asset heterogeneity and component heterogeneity, as well as economic dependencies between assets. For tractability purposes, the current formulation does not incorporate resource constraints. However, it can be extended to include shared resources\textemdash such as repair crews or spare parts\textemdash thereby allowing its application to asset networks with additional logistical complexities. Moreover, the degradation dynamics are modeled as a compound Poisson process with independent and identically distributed shock sizes. While this choice is standard in the reliability literature and supported by our case study data, it abstracts from certain practical complexities. In particular, we assume independent shock processes across assets and stationarity of degradation parameters within a component’s lifetime, thereby excluding, for example, common-cause or environment-driven correlations and non-stationary or regime-switching degradation behavior. Extending the framework to correlated or non-stationary degradation processes constitutes a promising direction for future research.

Using a partially observable Markov decision process (POMDP) framework, our approach leverages real-time degradation data to learn cost-effective maintenance strategies for asset networks with economic dependencies. We further addressed the computational challenges of solving POMDPs by employing a deep reinforcement learning (DRL) approach, enabling the derivation of near-optimal maintenance policies under parameter uncertainty. The DRL-based approximate policy iteration algorithm demonstrated its effectiveness by learning complex opportunistic maintenance strategies, directly improving upon common heuristic methods.

Through theoretical contributions, including the establishment of the structural properties of optimal replacement policies of the underlying Markov decision process (MDP) model and the reformulation of the POMDP as a Bayesian MDP (BMDP), our approach facilitates scalable and efficient maintenance solutions for industrial-scale asset networks. However, the scope of the established structural properties is limited to the full-information setting and does not directly transfer to the partial-information case. Nevertheless, we expect that the monotonicity result in the degradation level also holds under partial information, and a similar monotonicity result could be established for another state variable, namely the operational age, which we leave for future work. Such structural properties could aid the development of heuristic methods tailored to the BMDP case.

The practical value of the proposed model was highlighted through a case study on degrading interventional X-ray system components, providing actionable managerial insights. The instances considered in this work are relatively small and assume asset-level homogeneity (i.e., all component replacements stem from the same pool), and therefore require relatively few samples to train our DRL algorithm; industrial-scale problems require many more, which can make single-node execution potentially prohibitive. By leveraging distributed computing and high-performance clusters, sample collection can be parallelized, keeping computation times manageable even for large, complex industrial-scale instances. Our approach benefits significantly from smart initialization heuristics; developing a scalable initialization method for cases where heterogeneity is present at both the asset and component levels therefore appears promising.

However, when no computationally convenient parametric form exists (i.e., prior-to-posterior belief updating is not tractable) to model component heterogeneity, exact reformulation of the maintenance problem as a BMDP is no longer feasible. Alternative inference techniques, such as variational inference or Monte Carlo sampling methods, introduce approximations and significant computational challenges in the DRL algorithm. Consequently, it is not clear what the effect is on the trained DRL-based policies using such techniques. As a workaround, we proposed a solution where policies trained on the underlying MDP are applied in a POMDP setting using an open-loop feedback approach. However, our numerical experiments indicate that this solution is effective only in controlled settings and not when component heterogeneity is highly volatile, presenting an opportunity for future work.

Finally, the performance results of the trained policies under different information levels provide insight into the value of increasing the information level. In practice, additional information is obtained through investments in data availability, quality and processing infrastructure. The value of information can therefore be interpreted as a trade-off between such investment costs and the potential for improved CBM planning. Asset managers can evaluate this trade-off through scenario analysis to support investment decisions, which is a promising direction for future research.

\small
\paragraph{Acknowledgements} This work used the Dutch national e-infrastructure with the support of the SURF Cooperative using grant no. EINF-5192. This work was supported by the Netherlands Organisation for Scientific Research (NWO). Project: NWO Big data - Real Time ICT for Logistics. Number: 628.009.012. The work of Stella Kapodistria is supported by NWO through the Gravitation-grant NETWORKS-024.002.003. Collin Drent received support from NWO through Grant VI.Veni.241E.058.
\normalsize	
\bibliography{references.bib}
\pagebreak

\setcounter{page}{1}
\section*{Supplementary material}
\beginsupplement
\appendix
\setcounter{theorem}{0}
\section{Theorems and proofs}\label{app:proofs}

We provide the proofs of the established properties of the optimal replacement policy of the \emph{underlying} MDP model presented in \mbox{Section \ref{subsec:structural_properties}}, inspired by the proof of~\citet[Theorem 3.7]{soltani2023structured}. Note that we restrict attention to the case of $M=2$ parallel assets; however, the arguments extend analogously to the case $M>2$.

\begin{theorem}[Monotonicity in degradation levels]\label{thm:deglvls_app}
Let $x = (x_1, \lambda_1, \phi_1, \ldots, x_M, \lambda_M, \phi_M)$ denote a state. Define $\mathcal{M}(x) \subseteq \mathcal{M}$ to be the set of machines for which it is optimal to do maintenance in state $x$. For any $y = (y_1, \lambda_1, \phi_1, \ldots, y_M, \lambda_M, \phi_M)$ that represents a state with more severe degradation than $x$, i.e.,
\begin{enumerate}
 \item $y_m \geq x_m$ for all $m \in \mathcal{M}(x)$, and
 \item $y_m = x_m$ otherwise.
\end{enumerate}
Then, it holds that
\begin{equation*}
 a^*(x) \in \textrm{arg} \underset{a\in\mathcal{U}(y)}{\textrm{min}} Q(y,a).
\end{equation*}

Here, $a^*(x) \in \mathcal{U}(x)$ and $a^*(y) \in \mathcal{U}(y)$ denote the optimal actions in states $x$ and $y$, respectively, and $Q(\cdot,\cdot)$ denotes the state-action value function of an optimal policy. In other words, the optimal action $a^*(x)$ for state $x$ is also an optimal action for a state $y$ with more severe degradation.

\begin{proof}
Firstly, note that $\mathcal{U}(y) \subseteq \mathcal{U}(x)$ and that $a^*(x) \in \mathcal{U}(y)$. To prove the statement, we use a one-step argument to expand the LHS and show that, for all $a \in \mathcal{U}(y)$, the RHS is nonnegative:
\begin{align*}
Q(y,a) - Q(y, a^*(x)) &= \underbrace{C(y,a) - C(y, a^*(x))}_{\text{cost difference}} + \underbrace{\gamma\mathbb{E}[ V(\hat{y}(a)) ] - \gamma\mathbb{E}[ V(\hat{y}(a^*(x)))]}_{\text{future-value difference}}.
\end{align*}
We proceed with deriving bounds for each bracketed term on the RHS.

For $a \in \mathcal{U}(x)$, define the following sets:
\begin{equation*}
 \mathcal{J}(x,a) = \{m \in \mathcal{M} \mid a^*_m(x) = 1, a_m = 0 \},
\end{equation*}
\begin{equation*}
 \mathcal{G}(x,a) = \{m \in \mathcal{M} \mid a^*_m(x) = 0, a_m = 1 \} \textrm{ and}
\end{equation*}
By definition of $a^*(x)$,
\begin{equation*}
 Q(x, a^*(x)) \leq Q(x,a) \textrm{ for all } a\in\mathcal{U}(x).
\end{equation*}
We first show that for all $a\in\mathcal{U}(y)$
\begin{equation}\label{eq:app:thm1-1}
 C(y,a) - C(y, a^*(x)) \geq C(x,a) - C(x, a^*(x)).
\end{equation}
Corrective maintenance on any failed machine $m$ in state $y$ is mandatory, i.e., $\mathcal{U}_m(y) = \{1\}$. Therefore,
\begin{align*}
 C(y,a) - C(y, a^*(x)) &= c^\textrm{ST}\left( \mathds{1}_{\{\sum_{m\in\mathcal{M}} a_m > 0\}} - \mathds{1}_{\{ \sum_{m\in\mathcal{M}} a_m^*(x)>0 \}} \right) + \sum_{m\in \mathcal{G}(x,a)} c_m^\textrm{PM} - \sum_{m\in \mathcal{J}(x,a) } c_m^\textrm{PM} \\
 &= C(x,a) - C(x, a^*(x)).
\end{align*}
This implies \mbox{Eq. \ref{eq:app:thm1-1}}.

 Denote with $\hat{y}(a)$ and $\hat{x}(a)$ the next state after taking action $a \in \mathcal{U}(y)$ in states $y$ and $x$, respectively. Since we do maintenance on all machines with a more severe degradation level, we have that
\begin{equation*}
 \hat{y}(a^*(x)) \overset{\textrm{d}}{=} \hat{x}(a^*(x)).
\end{equation*}
For $a \in \mathcal{U}(y)$, note that $\hat{y}(a) \geq_\textrm{st} \hat{x}(a)$ element-wise. Equivalently,
\begin{equation*}
 \mathbb{E}[\phi(\hat{y}(a))] \geq \mathbb{E}[ \phi(\hat{x}(a))] 
\end{equation*}
for all \emph{non-decreasing} functions $\phi:\mathbb{R} \rightarrow \mathbb{R}$ for which these expectations exist~\citep[1.A.7]{shaked2007stochastic}. Note that it suffices to have this property for one state variable as the argument can be repeated for other state variables. In particular, we can show that this property holds for the value function.

Let $V(\cdot)$ denote the value function of an optimal policy that solves \mbox{Eq. (\ref{eq:objectiveMDP})} and let $V^n(\cdot)$ denote the value function at the $n$-th iteration of the value iteration algorithm, i.e., the value iteration algorithm produces the sequence $\{ V^n(\cdot) \}_{n\in\mathbb{Z}_+}$. Since our state space is a Borel space, $\gamma \in (0,1)$ and costs are bounded from above, the value iteration algorithm is guaranteed to converge point-wise to the optimal value function of an optimal policy~\citep[Proposition 9.14]{bertsekas1996stochastic}. That is, for all states $h$ of the underlying MDP, $V^n(h) \rightarrow V(h)$ for any arbitrary starting position through the value iteration algorithm.

For illustration purposes, we restrict ourselves to showing that the value function is non-decreasing with respect to the degradation level in the case of $M = 2$, although the argument can be extended to cases where $M>2$. We apply induction on the iterations of the value iteration algorithm to prove this for $x_1$. Note that it suffices to prove this for one degradation level as the argument can be repeated for the other degradation levels. For all states $h$, we initialize $V^0(h) = 0$. By definition, $V^0(\cdot)$ is non-decreasing in $x_1$. Assume that $V^n(\cdot)$, $n\in\mathbb{N}$, is non-decreasing in $x_1$.

To avoid notational clutter, we will omit the degradation process parameters $\lambda_1, \phi_1, \lambda_2$, and $\phi_2$ in the following argument. Let $D(t) = (X_1(t), X_2(t))$ denote the random amount of damage accumulated per machine after $t$ time units (without maintenance intervention). Let $s_0 = (0,0)$ be an arbitrary healthy state. The Bellman optimality equations for the case $M=2$ read as follows:\\
\begin{equation*}
\resizebox{\textwidth}{!}{$
V( (x_1, x_2) ) =
\begin{cases}
c^\textrm{ST} + c_1^\textrm{CM} + c_2^\textrm{CM} + \gamma \mathbb{E}\big[V( s_0 + D(1) )\big] & \text{if $(x_1, x_2) \geq (\xi_1, \xi_2)$,} \\
c^\textrm{ST} + c_1^\textrm{CM} + \textrm{min}\Big\{ \gamma \mathbb{E}\big[V( (0, x_2) + D(1) )\big], c_2^\textrm{PM} + \gamma \mathbb{E}\big[V( s_0 + D(1) )\big]\Big\} & \text{if $x_1 \geq \xi_1$ and $x_2 < \xi_2$,} \\
c^\textrm{ST} + c_2^\textrm{CM} + \textrm{min}\Big\{ \gamma \mathbb{E}\big[V( (x_1, 0) + D(1) )\big], c_1^\textrm{PM} + \gamma \mathbb{E}\big[V( s_0 + D(1) )\big]\Big\} & \text{if $x_1 < \xi_1$ and $x_2 \geq \xi_2$,} \\
\textrm{min} \Big\{ \gamma \mathbb{E}\big[V( (x_1, x_2) + D(1) )\big], c^\textrm{ST} + c_2^\textrm{PM} + \gamma \mathbb{E}\big[V( (x_1, 0) + D(1) )\big], \\ c^\textrm{ST} + c_1^\textrm{PM} + \gamma \mathbb{E}\big[V( (0, x_2) + D(1) )\big], c^\textrm{ST} + c_1^\textrm{PM} +c_2^\textrm{PM} + \gamma \mathbb{E}\big[V( s_0 + D(1) )\big]\Big\} & \text{if $(x_1, x_2) < (\xi_1, \xi_2)$.} \end{cases}
$}
\end{equation*}

From the Bellman equations, we have that
\begin{equation*}
\resizebox{\textwidth}{!}{$
V^{n+1}( (x_1, x_2) ) =
\begin{cases}
c^\textrm{ST} + c_1^\textrm{CM} + c_2^\textrm{CM} + \gamma \mathbb{E}\big[V^n( s_0 + D(1) )\big] & \text{if $(x_1, x_2) \geq (\xi_1, \xi_2)$,} \\
c^\textrm{ST} + c_1^\textrm{CM} + \textrm{min}\Big\{ \gamma \mathbb{E}\big[V^n( (0, x_2) + D(1) )\big], c_2^\textrm{PM} + \gamma \mathbb{E}\big[V^n( s_0 + D(1) )\big]\Big\} & \text{if $x_1 \geq \xi_1$ and $x_2 < \xi_2$,} \\
c^\textrm{ST} + c_2^\textrm{CM} + \textrm{min}\Big\{ \gamma \mathbb{E}\big[V^n( (x_1, 0) + D(1) )\big], c_1^\textrm{PM} + \gamma \mathbb{E}\big[V^n( s_0 + D(1) )\big]\Big\} & \text{if $x_1 < \xi_1$ and $x_2 \geq \xi_2$,} \\
\textrm{min} \Big\{ \gamma \mathbb{E}\big[V^n( (x_1, x_2) + D(1) )\big], c^\textrm{ST} + c_1^\textrm{PM} + \gamma \mathbb{E}\big[V^n( (0, x_2) + D(1) )\big], \\c^\textrm{ST} + c_2^\textrm{PM} + \gamma \mathbb{E}\big[V^n( (x_1, 0) + D(1) )\big], c^\textrm{ST} + c_1^\textrm{PM} +c_2^\textrm{PM} + \gamma \mathbb{E}\big[V^n( s_0 + D(1) )\big]\Big\} & \text{if $(x_1, x_2) < (\xi_1, \xi_2)$.} \end{cases}
$}
\end{equation*}
All the constants that appear on the RHS and $D(1)$ are non-decreasing in $x_1$; therefore, $V^{n+1}( (x_1, x_2) )$ is non-decreasing in $x_1$. All in all, this implies that
\begin{equation}\label{eq:app:thm1-2}
 \mathbb{E}[V(\hat{y}(a))] \geq \mathbb{E}[V(\hat{x}(a))].
\end{equation}
We proceed by proving that for all $a\in\mathcal{U}(y)$,
\begin{equation*}
 Q(y,a) - Q(y, a^*(x)) \geq Q(x,a) - Q(x, a^*(x)).
\end{equation*}
This result follows directly from the previously established results in \mbox{Eq. (\ref{eq:app:thm1-1})} and \mbox{Eq. (\ref{eq:app:thm1-2})}:
\begin{align*}
 Q(y,a) - Q(y, a^*(x)) &= C(y,a) - C(y, a^*(x)) + \gamma\mathbb{E}[ V(\hat{y}(a)) ] - \gamma\mathbb{E}[ V(\hat{y}(a^*(x))) ] \\
 &\geq C(x,a) - C(x, a^*(x)) + \gamma\mathbb{E}[ V(\hat{x}(a)) ] - \gamma\mathbb{E}[ V(\hat{x}(a^*(x))) ] \\
 &= Q(x,a) - Q(x, a^*(x)) \geq 0.
\end{align*}
This implies that
\begin{equation*}
 Q(y,a) \geq Q(y,a^*(x)).
\end{equation*}
In other words,
\begin{equation*}
 a^*(x) \in \textrm{arg} \underset{a\in\mathcal{U}(y)}{\textrm{min}} Q(y,a).
\end{equation*}
By definition,
\begin{equation*}
 a^*(y) \in \textrm{arg} \underset{a\in\mathcal{U}(y)}{\textrm{min}} Q(y,a).
\end{equation*}
This means that
\begin{equation*}
 Q(y,a^*(x)) = Q(y,a^*(y)).
\end{equation*}
This concludes the proof.
\end{proof}
\end{theorem}

\begin{theorem}[Monotonicity in degradation parameters]\label{thm:degparams_app}
Let $x = (x_1, \lambda_1, \phi_1, \ldots, x_M, \lambda_M, \phi_M)$ denote a state. Let $\mathcal{M}(x) \subseteq \mathcal{M}$ denote the set of machines for which it is optimal to do maintenance in state $x$. For any $x^\prime = (x_1, \lambda^\prime_1, \phi^\prime_1, \ldots, x_M, \lambda^\prime_M, \phi^\prime_M)$ that represents a state with worse degradation parameters than $x$, i.e.,
\begin{enumerate}
 \item $\lambda^\prime_m \geq \lambda_m$ and $\phi^\prime_m \geq \phi_m$ for all $m \in \mathcal{M}(x)$, and
 \item $\lambda^\prime_m = \lambda_m$ and $\phi^\prime_m = \phi_m$ otherwise.
\end{enumerate}

If $\phi'_m \geq \phi_m$ implies that the corresponding shock size distributions satisfy $Y'_m \geq_\textrm{st} Y_m$ (where $Y'_m$ and $Y_m$ are the random variables associated with the respective shock size distributions), then it holds that
\begin{equation*}
 a^*(x) \in \textrm{arg} \underset{a\in\mathcal{U}(x^\prime)}{\textrm{min}} Q(x^\prime,a).
\end{equation*}

The optimal action $a^*(x)$ for state $x$ is thus also an optimal action for state $x^\prime$.

\begin{proof}
Note that, in this case, $\mathcal{U}(x^\prime) = \mathcal{U}(x)$ and thus $a^*(x) \in \mathcal{U}(x^\prime)$. The proof follows the same approach as the proof of \mbox{Theorem \ref{thm:deglvls_app}}. Using a one-step argument, we expand the LHS and show that, for all $a \in \mathcal{U}(x^\prime)$, the RHS is nonnegative:
\begin{align*}
Q(x^\prime,a) - Q(x^\prime, a^*(x)) &= \underbrace{C(x^\prime,a) - C(x^\prime, a^*(x))}_{\text{cost difference}} + \underbrace{\gamma\mathbb{E}[ V(\hat{x}^\prime(a)) ] - \gamma\mathbb{E}[ V(\hat{x}^\prime(a^*(x)))]}_{\text{future-value difference}}.
\end{align*}
We proceed with deriving bounds for each bracketed term on the RHS.

Since states $x$ and $x^\prime$ have identical degradation levels and differ only in their degradation process parameters, the costs of taking any action are the same. Therefore, for all $a \in \mathcal{U}(x)$,
\begin{equation*}
 C(x^\prime,a) = C(x,a),
\end{equation*}
which implies, in particular,
\begin{equation*}
 C(x^\prime, a^*(x)) = C(x, a^*(x)).
\end{equation*}
Thus, we have
\begin{equation}\label{eq:app:thm2-1}
 C(x^\prime,a) - C(x^\prime, a^*(x)) = C(x,a) - C(x, a^*(x)).
\end{equation}
Denote with $\hat{x}^\prime(a)$ and $\hat{x}(a)$ the next state after taking action $a \in \mathcal{U}(x)$ in states $x^\prime$ and $x$, respectively. Since we do maintenance on all machines with ``worse'' degradation process parameters, we have that
\begin{equation*}
 \hat{x}^\prime(a^*(x)) \overset{\textrm{d}}{=} \hat{x}(a^*(x)).
\end{equation*}
The assumptions $\lambda^{\prime}_m \geq \lambda_m$ and $Y^{\prime}_m \geq_\textrm{st} Y_m$ for all $m \in \mathcal{M}(x)$ ensure that $\hat{x}^\prime(a) \geq_\textrm{st} \hat{x}(a)$ for all $a \in \mathcal{U}(x)$~\citep[Theorem 1.A.4]{shaked2007stochastic}. A similar argument to the one used in the proof of \mbox{Theorem \ref{thm:deglvls}}, which uses induction on the iterations of the value iteration algorithm, shows that the value function is non-decreasing with respect to the degradation process parameters $\lambda_m$ and $\phi_m$. This means that
\begin{equation}\label{eq:app:thm2-2}
 \mathbb{E}[V(\hat{x}^\prime(a))] \geq \mathbb{E}[V(\hat{x}(a))]. 
\end{equation}
The remainder of the proof is identical to the proof of \mbox{Theorem \ref{thm:deglvls}}. We show that
\begin{equation*}
 Q(x^\prime,a) - Q(x^\prime, a^*(x)) \geq Q(x,a) - Q(x, a^*(x)).
\end{equation*}
This result follows directly from the previously established results in \mbox{Eq. (\ref{eq:app:thm2-1})} and \mbox{Eq. (\ref{eq:app:thm2-2})}:
\begin{align*}
 Q(x^\prime,a) - Q(x^\prime, a^*(x)) &= C(x^\prime,a) - C(x^\prime, a^*(x)) + \gamma\mathbb{E}[ V(\hat{x}^\prime(a)) ] - \gamma\mathbb{E}[ V(\hat{x}^\prime(a^*(x))) ] \\
 &\geq C(x,a) - C(x, a^*(x)) + \gamma\mathbb{E}[ V(\hat{x}(a)) ] - \gamma\mathbb{E}[ V(\hat{x}^\prime(a^*(x))) ] \\
  &= C(x,a) - C(x, a^*(x)) + \gamma\mathbb{E}[ V(\hat{x}(a)) ] - \gamma\mathbb{E}[ V(\hat{x}(a^*(x))) ] \\
 &= Q(x,a) - Q(x, a^*(x)) \geq 0.
\end{align*}
Therefore,
\begin{equation*}
 Q(x^\prime,a) \geq Q(x^\prime,a^*(x))
\end{equation*}
for all $a\in\mathcal{U}(x^\prime)$. In other words,
\begin{equation*}
 a^*(x) \in \textrm{arg} \underset{a\in\mathcal{U}(x^\prime)}{\textrm{min}} Q(x^\prime,a).
\end{equation*}
By definition,
\begin{equation*}
 a^*(x^\prime) \in \textrm{arg} \underset{a\in\mathcal{U}(x^\prime)}{\textrm{min}} Q(x^\prime,a).
\end{equation*}
This means that
\begin{equation*}
 Q(x^\prime,a^*(x)) = Q(x^\prime,a^*(x^\prime)).
\end{equation*}
This concludes the proof.
\end{proof}
\end{theorem}

\newpage

\section{Nomenclature} \label{app:nomenc}
\begin{table*}[!ht]
 \begin{framed}
 \printnomenclature
\nomenclature[01]{$\mathcal{M}$}{Set of machines}
\nomenclature[02]{$M$}{Number of machines}
\nomenclature[03]{$m$}{Machine index}
\nomenclature[04]{$t$}{Time step}
\nomenclature[05]{$N_m(t)$}{Shock arrival process of machine $m$}
\nomenclature[06]{$\lambda_m$}{Shock arrival rate of machine $m$}
\nomenclature[07]{$\Lambda_m$}{Sampling distribution of the shock \\arrival rate of machine $m$}
\nomenclature[08]{$\alpha_m$, $\beta_m$}{Shape and scale hyperparameters of the Gamma prior for $\lambda_m$}
\nomenclature[09]{$Y^{(i)}_m$}{Shock size of machine $m$}
\nomenclature[10]{$\phi_m$}{Shock size distribution parameter of machine $m$}
\nomenclature[11]{$p_m$}{Geometric distribution parameter of machine $m$}
\nomenclature[12]{$\Phi_m$}{Sampling distribution of the shock size distribution parameter of machine $m$}
\nomenclature[13]{$r_m$, $s_m$}{Hyperparameters of the conjugate prior for $\phi_m$}
\nomenclature[14]{$X_m(t)$}{Degradation process of machine $m$}
\nomenclature[15]{$x_m(t)$}{Degradation level of machine $m$}
\nomenclature[16]{$\xi_m$}{Failure threshold of machine $m$}
\nomenclature[17]{$k_m(t)$}{Shocks sustained by machine $m$ in \\period $(t-1,t]$}
\nomenclature[18]{$z_m(t)$}{Accumulated damage of machine $m$ in period $(t-1,t]$}
\nomenclature[19]{$t_m(t)$}{Machine age of machine $m$}
\nomenclature[20]{$h$}{Underlying MDP state}
\nomenclature[21]{$o$}{Observed state}
\nomenclature[22]{$\tilde{h}$}{BMDP state}
\nomenclature[23]{$\mathcal{U}(h)$}{State-dependent action set}
\nomenclature[24]{$c^\textrm{PM}_m$}{PM cost of machine $m$}
\nomenclature[25]{$c^\textrm{CM}_m$}{CM cost of machine $m$}
\nomenclature[26]{$c^\textrm{ST}$}{Setup cost}
\nomenclature[27]{$\gamma$}{Discount factor}
\nomenclature[28]{$C(h,a)$}{Cost function}
\nomenclature[29]{$\infolevel{i}$}{Information level $i$}
\nomenclature[30]{$\pi^{\infolevel{i}}$}{Policy under information level $\infolevel{i}$}
\nomenclature[31]{$J(\pi^{\infolevel{i}})$}{Total expected discounted cost under policy $\pi^{\infolevel{i}}$}
\nomenclature[32]{$\tau_m^\textrm{PM}$}{PM threshold of machine $m$}
\nomenclature[33]{$\tau_m^\textrm{OPM}$}{OPM threshold of machine $m$}
\nomenclature[34]{$f^{\infolevel{}}(h)$}{Feature vector}
 \end{framed}
\end{table*}

\newpage

\section{Examples of exponential-family representations}
\label{app:examples}

In this appendix, we present two commonly used discrete shock-size distributions and show how they can be written in canonical exponential-family form with a linear sufficient statistic.

\subsection{Binomial distribution}

Suppose the shock size is binomially distributed with a fixed number of trials $n_m \in \mathbb{N}$ and success probability $p_m \in (0,1)$. Its probability mass function is
\[
f(x \mid p_m)
= \binom{n_m}{x} p_m^x (1-p_m)^{n_m-x}, 
\quad x=0,1,\dots,n_m.
\]

Rewriting the density,
\[
f(x \mid p_m)
= \binom{n_m}{x}
\exp\Big( x \ln p_m + (n_m-x)\ln(1-p_m) \Big)
\]
\[
= \binom{n_m}{x}
\exp\Big( x \ln\!\Big(\frac{p_m}{1-p_m}\Big)
+ n_m \ln(1-p_m) \Big).
\]

Hence, the binomial distribution can be written in canonical exponential-family form
\[
f(x \mid \phi_m)
= h_m(x)\exp\big( \phi_m T_m(x) - A_m(\phi_m) \big),
\]
with
\[
T_m(x)=x,
\quad
\phi_m=\ln\!\Big(\frac{p_m}{1-p_m}\Big),
\]
\[
h_m(x)=\binom{n_m}{x},
\quad
A_m(\phi_m)=n_m \ln\big(1+e^{\phi_m}\big).
\]

Thus, the binomial distribution belongs to the one-parameter exponential family with linear sufficient statistic. 

\subsection{Poisson distribution}

Suppose the shock size is Poisson distributed with mean $\mu_m>0$. Its probability mass function is
\[
f(x \mid \mu_m)
= \frac{e^{-\mu_m}\mu_m^x}{x!},
\quad x=0,1,2,\dots.
\]

Rewriting,
\[
f(x \mid \mu_m)
= \frac{1}{x!}
\exp\big( x \ln \mu_m - \mu_m \big).
\]

This can be written in canonical exponential-family form as
\[
f(x \mid \phi_m)
= h_m(x)\exp\big( \phi_m T_m(x) - A_m(\phi_m) \big),
\]
with
\[
T_m(x)=x,
\quad
\phi_m=\ln \mu_m,
\]
\[
h_m(x)=\frac{1}{x!},
\quad
A_m(\phi_m)=e^{\phi_m}.
\]

Hence, the Poisson distribution also belongs to the one-parameter exponential family with linear sufficient statistic. 

\newpage
\section{Deep controlled learning hyperparameters}\label{app:drl_param}
We present the selected hyperparameters for the deep controlled learning algorithm below. All neural networks employ the rectified linear unit (ReLU) activation function.

\begin{table}[!ht]
\centering
\begin{tabular}{c|l|c}
\toprule
Hyperparameter & Description & Value \\ \midrule
$L$ & number of neural network layers & $4$ \\
$d^1$ & dimension of input layer & $256$ \\
$d^l$ & dimension of hidden layer $l=2, \ldots, L$ & $128$ \\
$\text{MAX\_SAMPLES}$ & number of samples & 50,000 (CS.3: 250,000) \\
$r_\textrm{max}$ & maximum number of roll-outs & $7500$ \\
$\text{BATCH\_SIZE}$ & batch size & $64$\\
$k$ & $k$-value bandit optimizer & $2.0$ \\
$\epsilon$ & fraction random actions & $0.02$ \\
 \bottomrule
\end{tabular}%
\caption{Approximate policy iteration hyperparameters for the considered instances.}
\label{tab: hyper_instances}
\end{table}

\clearpage
\section{Deep controlled learning training algorithm, accuracy and costs}\label{app:drl_training}

Deep controlled learning relies on supervised learning to train the neural network. In the $\infolevel{1}$ setting, i.e., the Bayesian Markov decision process formulation described in \mbox{Section \ref{subsec:bmdp}}, the goal is to learn the mapping from the feature representation $f^{\infolevel{1}}_3(\tilde{h}^{a_{\sigma(m-1)}})$ of a state $\tilde{h}^{a_{\sigma(m-1)}}$ to the action $\tilde{a}_* \in \mathcal{U}_{\sigma(m)}(\tilde{h}^{a_{\sigma(m-1)}})$ stored in the data set $\mathcal{D}$. This mapping is approximated using a multilayer perceptron with $L \in \mathbb{Z}_+$ layers, where each layer $l \in \{1,\ldots,L\}$ applies an affine transformation followed by a nonlinear activation function.

We adopt a standard supervised learning procedure; see~\citet[Appendix A]{temizoz2023deep}. The data set $\mathcal{D}$ is split into a training set and a validation set. The parameters $\theta$ of the neural network are estimated by minimizing a cross-entropy loss $L(\theta)$ that measures the discrepancy between the target policy and the neural network policy $\pi^{\infolevel{1}}_\theta$ on the training data. The optimization is performed iteratively using gradient-based updates: in each iteration, the gradient of $L(\theta)$ with respect to $\theta$ is computed and the parameters are updated in the direction of decreasing loss. Training is stopped once the validation loss, defined analogously, no longer improves.

\mbox{Table \ref{tab:trainingresults}} reports the classification accuracy obtained during training for the case study instances CS.1--CS.3. Since these instances are still relatively small, the performance is overall very high. We expect that, as problem complexity increases, representing a near-optimal policy via function approximation will become more challenging, which will likely reduce the classification accuracy of the neural network.

\begin{table}[H]
\centering
\begin{tabular}{c|c|c|c}
\toprule
 & CS.1 & CS.2 & CS.3 \\ \midrule
Validation set (gen 1): & 93\% & 95\% & 93\%   \\ 
Training set (gen 1): & 96\% & 96\% & 95\% \\ \hline
Validation set (gen 2): & 95\% & 94\%  & 93\% \\ 
Training set (gen 2): & 96\% & 95\% & 96\% \\ \hline
Validation set (gen 3): & 94\% & 94\% & 94\% \\ 
Training set (gen 3): & 96\% & 96\% & 95\%
\end{tabular}
\caption{Classification accuracy of the trained neural networks on the training and validation set for the case study instances CS.1--3.}
\label{tab:trainingresults}
\end{table}

\begin{figure}[!ht]
 \centering
 \begin{subfigure}[t]{0.33\textwidth}
 \centering
 \includegraphics[width=\textwidth]{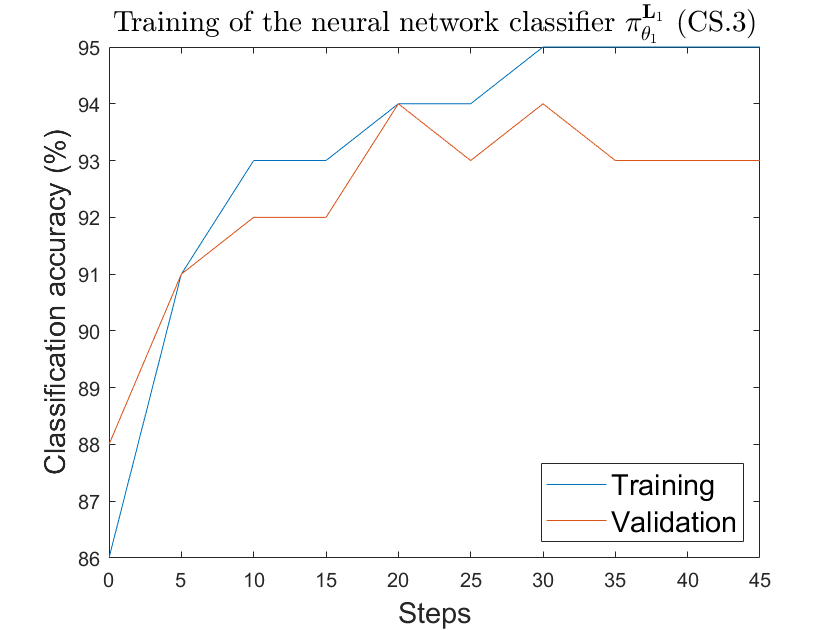}
 \label{subfig:CS3BMDP_classifier_gen1}
 \end{subfigure}%
 \hfill
 \begin{subfigure}[t]{0.33\textwidth}
 \centering
 \includegraphics[width=\textwidth]{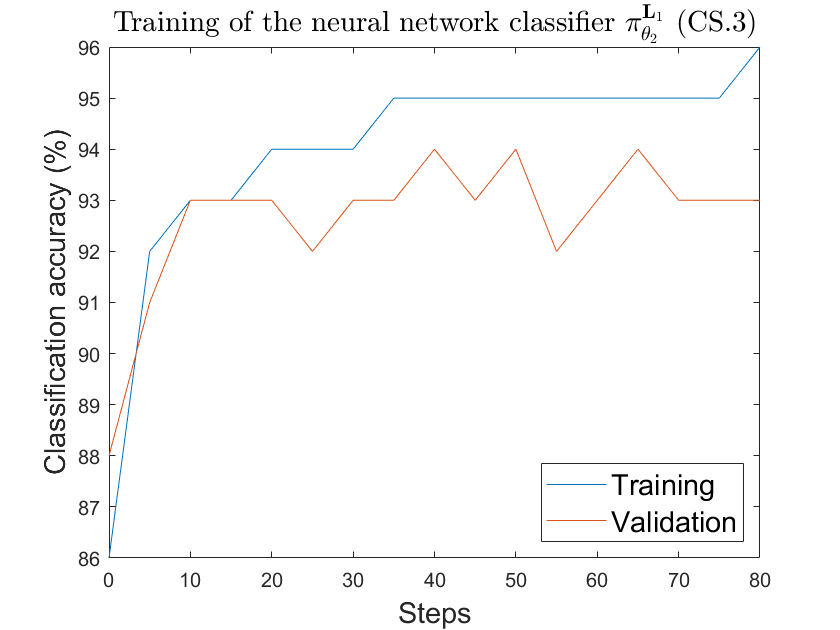}
 \label{subfig:CS3BMDP_classifier_gen2}
 \end{subfigure}
 \hfill
 \begin{subfigure}[t]{0.33\textwidth}
 \centering
 \includegraphics[width=\textwidth]{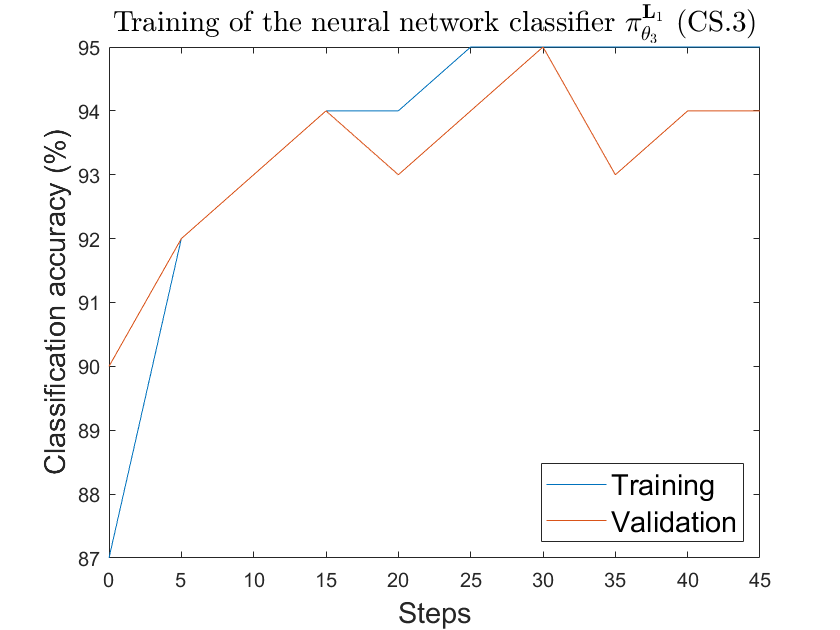}
 \label{subfig:CS3BMDP_classifier_gen3}
 \end{subfigure}
\caption{Accuracy of the neural networks trained in the $\infolevel{1}$ setting for the case study instance CS.3.}
 \label{fig:api_classifier}
\end{figure}

We illustrate in \mbox{Figure \ref{fig:api_classifier}} the accuracy of the neural network on the training and validation sets during training. The neural network rapidly achieves high classification accuracy, typically within the first 20 steps. The accuracy on the training data reaches 95--96\%, while accuracy on the validation set stabilizes around 92--94\%, indicating good performance on unseen data.

Finally, \mbox{Table \ref{tab:trainingtimes}} reports the cost and duration of the sample collection phase for case study instances CS.1–3. These instances have relatively small state spaces, which in turn limits the number of required samples. However, for industrial-scale problem instances, a substantially greater number of samples may be necessary to adequately explore the state space. This increase in sample size can impose a significant computational burden, and execution of the algorithm on a single node may thus become prohibitively slow. To address this, distributed computing can be leveraged, enabling the sample collection phase to be parallelized across multiple nodes. By utilizing high-performance computing clusters, it is possible to maintain manageable sample collection times even for large and complex instances.
\begin{table}[H]
\centering
\begin{tabular}{c|c|c|c}
\toprule
 & CS.1 & CS.2 & CS.3 \\ \midrule
Training time (gen 1) & 0.076 hrs & 0.140 hrs & 1.661 hrs  \\
Cost & \euro0.22 & \euro0.40 & \euro4.78 \\ \midrule
Training time (gen 2) & 0.196 hrs & 0.412 hrs & 5.135 hrs \\ 
Cost & \euro0.56 & \euro1.19 & \euro14.79 \\ \midrule
Training time (gen 3) & 0.196 hrs & 0.408 hrs & 5.114 hrs \\ 
Cost & \euro0.56 & \euro1.18 & \euro14.73  
\end{tabular}
\caption{Cost and duration of the sample collection phase. Each reported time represents an estimated duration required to collect the necessary number of samples (see \ref{app:drl_param}) using a single \emph{thin computing node} on the Dutch National Supercomputer~\citep{snellius}. The associated cost is calculated based on~\citet{surf2024}.}
\label{tab:trainingtimes}
\end{table}

\end{document}